# GENERALIZED WITTEN GENUS AND VANISHING THEOREMS

QINGTAO CHEN, FEI HAN, AND WEIPING ZHANG

ABSTRACT. We construct a generalized Witten genus for spin$^c$ manifolds, which takes values in level 1 modular forms with integral Fourier expansion on a class of spin$^c$ manifolds called string$^c$ manifolds. We also construct a mod 2 analogue of the Witten genus for $8k+2$ dimensional spin manifolds. The Landweber-Stong type vanishing theorems are proven for the generalized Witten genus and the mod 2 Witten genus on string$^c$ and string (generalized) complete intersections in (product of) complex projective spaces respectively.

## 1. INTRODUCTION

Let $M$ be a $4k$ dimensional closed oriented smooth manifold. Let $\{\pm 2\pi\sqrt{-1}z_j, 1 \leq j \leq 2k\}$ denote the formal Chern roots of $T_{\mathbf{C}}M$, the complexification of the tangent vector bundle $TM$ of $M$. Then the famous Witten genus of $M$ can be written as (cf. [19])

$$(1.1) \qquad W(M) = \left\langle \left(\prod_{j=1}^{2k} z_j \frac{\theta'(0,\tau)}{\theta(z_j,\tau)}\right), [M] \right\rangle \in \mathbf{Q}[[q]],$$

with $\tau \in \mathbf{H}$, the upper half-plane, and $q = e^{\pi\sqrt{-1}\tau}$.

The above genus was first introduced in [22] and can be viewed as the loop space analogue of the $\widehat{A}$-genus. It can be expressed as a $q$-deformed $\widehat{A}$-genus as

$$W(M) = \left\langle \widehat{A}(TM) \mathrm{ch}\left(\Theta\left(T_{\mathbf{C}}M\right)\right), [M] \right\rangle,$$

where

$$\Theta(T_{\mathbf{C}}M) = \overset{+\infty}{\underset{m=1}{\otimes}} S_{q^{2m}}(\widetilde{T_{\mathbf{C}}M}), \quad \text{with } \widetilde{T_{\mathbf{C}}M} = T_{\mathbf{C}}M - \mathbf{C}^{4k},$$

is the Witten bundle defined in [22].

The manifold $M$ is called spin, if $w_1(TM) = 0$, $w_2(TM) = 0$, where $w_1(TW)$, $w_2(TM)$ are the first and the second Stiefel-Whitney classes of $TM$ respectively. According to the Atiyah-Singer index theorem [3], if $M$ is spin, then $W(M) = \mathrm{Ind}(D \otimes \Theta(T_{\mathbf{C}}M)) \in \mathbf{Z}[[q]]$, where $D$ is the Atiyah-Singer-Dirac operator on $M$ (cf. [14]).

Moreover, $M$ is called string if one further topological restriction is put on $M$, i.e. $\frac{1}{2}p_1(TM) = 0$, where $p_1(TM)$ is the first Pontrjagin class of $TM$.[1] It is well-known that if $M$ is string, or even weaker, if $M$ is spin and the first rational Pontryagin class of $M$ vanishes, then $W(M)$ is a modular form of weight $2k$ over $SL(2,\mathbf{Z})$ with integral Fourier expansion ([23]). The homotopy theoretical refinements of the Witten genus on string manifolds leads to the theory of tmf (topological modular form) developed by Hopkins and Miller [15].

The *Witten operator* $D \otimes \Theta(T_{\mathbf{C}}M)$ was proved to be rigid, in the sense of Witten [22], by Liu in [18]. Moreover, Liu showed that if $M$ admits an $S^1$-action and $p_1(M)_{S^1} = n \cdot \pi^* u^2$, where $p_1(M)_{S^1}$ is the equivariant first Pontrjagin class, $\pi : M \times_{S^1} ES^1 \to BS^1$ is the projection, $u \in H^2(BS^1, \mathbf{Z})$ is the generator and $n$ is an integer, then the index of the Witten operator is

---
*Date*: March 11, 2010.
[1]The class $\frac{p_1(TM)}{2}$ is a degree 4 cohomology class, twice of which equals to $p_1(TM)$.

zero. We recall that Liu used in [18] the modularity of the Lefschetz numbers of certain twisted elliptic operators to prove the rigidity of the Witten operator as well as many other operators.

On the other hand, the Witten genus has a well-known vanishing result due to Landweber and Stong stating that $W(X) = 0$ if $X$ is a string complete intersection in a complex projective space (cf. [14, pp. 87-88]). In [7], Chen and Han generalize the result of Landweber and Stong by proving that the Witten genus of a string generalized complete intersection in a product of complex projective spaces vanishes.[2]

In this paper, we extend the classical Witten genus in two directions.

In one direction, we construct a Witten type genus for even dimensional spin$^c$ manifolds, which takes values in level 1 modular forms over $SL(2, \mathbf{Z})$ with integral Fourier expansions when the spin$^c$ manifolds are put more topological condition to become *string$^c$ manifolds*, just like the classical Witten genus takes values in level 1 modular forms over $SL(2, \mathbf{Z})$ with integral Fourier expansions on string manifolds. Similar to the classical Witten genus, which is the index of the Witten operator, this generalized Witten genus for spin$^c$ manifolds is the index of a Witten type operator, a twisted spin$^c$ Dirac operator. This Witten type operator can be shown to be rigid under relevant anomaly cancellation condition by applying Liu's method in [18]. We prove Landweber-Stong type vanishing theorems for this generalized Witten genus on string$^c$ (generalized) complete intersections in (product of) complex projective spaces, extending the vanishing theorems of Landweber-Stong and Chen-Han [7] on the classical Witten genus. It would be interesting to study the homotopy theoretical refinements of this generalized Witten genus for string$^c$ manifolds.

In the other direction, for $8k+2$-dimensional spin manifolds, we define a mod 2 analogue of the classical Witten genus. It is the mod 2 index of the Witten operator on $8k+2$-dimensional spin manifolds. We also obtain mod 2 analogues of the Landweber-Stong vanishing theorem of the classcial Witten genus on string complete intersections by proving vanishing theorems for this mod 2 Witten genus on $8k+2$-dimensional string (generalized) complete intersections in (product of) complex projective spaces.

The rest of the paper is organized as follows. In Section 2, we recall some knowledge of the Jacobi theta functions, modular forms as well as the definitions of characteristic forms to be discussed. In Section 3, we construct and discuss the generalized Witten genus for spin$^c$ manifolds as well as the mod 2 Witten genus for $8k+2$-dimensional spin manifolds. We then present and prove the Landweber-Stong type vanishing theorems for string$^c$ and string (generalized) complete intersections in (product of) complex projective spaces in Section 4.

Some of the results of this paper have been announced in [8].

## 2. Preliminaries

2.1. **Modular forms and Jacobi theta functions.** In this subsection, we recall some necessary knowledge on theta-functions and modular forms.

The four Jacobi theta-functions are defined as follows (cf. [5]),

$$(2.1) \qquad \theta(z,\tau) = 2q^{1/4}\sin(\pi z)\prod_{j=1}^{\infty}[(1-q^{2j})(1-e^{2\pi\sqrt{-1}z}q^{2j})(1-e^{-2\pi\sqrt{-1}z}q^{2j})],$$

$$(2.2) \qquad \theta_1(z,\tau) = 2q^{1/4}\cos(\pi z)\prod_{j=1}^{\infty}[(1-q^{2j})(1+e^{2\pi\sqrt{-1}z}q^{2j})(1+e^{-2\pi\sqrt{-1}z}q^{2j})],$$

---

[2]Recall that a generalized complete intersection in a compact oriented even-dimensional manifold $M$ is the transversal intersection of codimension two oriented submanifolds in $M$.

(2.3) $$\theta_2(z,\tau) = \prod_{j=1}^{\infty}[(1-q^{2j})(1-e^{2\pi\sqrt{-1}z}q^{2j-1})(1-e^{-2\pi\sqrt{-1}z}q^{2j-1})],$$

(2.4) $$\theta_3(z,\tau) = \prod_{j=1}^{\infty}\left[(1-q^{2j})(1+e^{2\pi\sqrt{-1}z}q^{2j-1})(1+e^{-2\pi\sqrt{-1}z}q^{2j-1})\right],$$

where $q = e^{\pi\sqrt{-1}\tau}$, $\tau \in \mathbf{H}$, the upper half plane.

They are all holomorphic functions for $(z,\tau) \in \mathbf{C} \times \mathbf{H}$, where $\mathbf{C}$ is the complex plane.

Let $\theta'(0,\tau) = \frac{\partial}{\partial z}\theta(z,\tau)|_{z=0}$, then one has the following Jacobi identity.

**Proposition 2.1** (Jacobi identity, [5])**.** *The following identity holds,*

(2.5) $$\theta'(0,\tau) = \pi\theta_1(0,\tau)\theta_2(0,\tau)\theta_3(0,\tau).$$

Let
$$SL(2,\mathbf{Z}) := \left\{ \begin{pmatrix} a_1 & a_2 \\ a_3 & a_4 \end{pmatrix} \middle| a_1, a_2, a_3, a_4 \in \mathbf{Z},\ a_1 a_4 - a_2 a_3 = 1 \right\}$$
be the modular group. Let $S = \begin{pmatrix} 0 & -1 \\ 1 & 0 \end{pmatrix}$, $T = \begin{pmatrix} 1 & 1 \\ 0 & 1 \end{pmatrix}$ be the two generators of $SL(2,\mathbf{Z})$. Their actions on $\mathbf{H}$ are given by

$$S : \tau \mapsto -\frac{1}{\tau}, \quad T : \tau \mapsto \tau + 1.$$

If we act theta-functions by $S$ and $T$, the following transformation formulas hold (cf. [5]),

(2.6) $$\theta(z,\tau+1) = e^{\frac{\pi\sqrt{-1}}{4}}\theta(z,\tau), \quad \theta(z,-1/\tau) = \frac{1}{\sqrt{-1}}\left(\frac{\tau}{\sqrt{-1}}\right)^{1/2} e^{\pi\sqrt{-1}\tau z^2}\theta(\tau z,\tau) ;$$

(2.7) $$\theta_1(z,\tau+1) = e^{\frac{\pi\sqrt{-1}}{4}}\theta_1(z,\tau), \quad \theta_1(z,-1/\tau) = \left(\frac{\tau}{\sqrt{-1}}\right)^{1/2} e^{\pi\sqrt{-1}\tau z^2}\theta_2(\tau z,\tau) ;$$

(2.8) $$\theta_2(z,\tau+1) = \theta_3(z,\tau), \quad \theta_2(z,-1/\tau) = \left(\frac{\tau}{\sqrt{-1}}\right)^{1/2} e^{\pi\sqrt{-1}\tau z^2}\theta_1(\tau z,\tau) ;$$

(2.9) $$\theta_3(z,\tau+1) = \theta_2(z,\tau), \quad \theta_3(z,-1/\tau) = \left(\frac{\tau}{\sqrt{-1}}\right)^{1/2} e^{\pi\sqrt{-1}\tau z^2}\theta_3(\tau z,\tau) .$$

There are also the following formulas,

(2.10) $$\theta(z+1,\tau) = -\theta(z,\tau),\ \theta(z+\tau,\tau) = -\frac{1}{q}e^{-2\pi i z}\theta(z,\tau),$$

(2.11) $$\theta_1(z+1,\tau) = -\theta_1(z,\tau),\ \theta_1(z+\tau,\tau) = \frac{1}{q}e^{-2\pi i z}\theta_1(z,\tau),$$

(2.12) $$\theta_2(z+1,\tau) = \theta_2(z,\tau),\ \theta_2(z+\tau,\tau) = -\frac{1}{q}e^{-2\pi i z}\theta_2(z,\tau),$$

(2.13) $$\theta_3(z+1,\tau) = \theta_3(z,\tau),\ \theta_3(z+\tau,\tau) = \frac{1}{q}e^{-2\pi i z}\theta_3(z,\tau).$$

Therefore it's not hard to deduce how the theta functions vary along the lattice $\Gamma = \{a + b\tau | a, b \in \mathbf{Z}\}$. In fact, we have

(2.14) $$\theta(z+a,\tau) = (-1)^a \theta(z,\tau)$$

and

$$
\begin{aligned}
&\theta(z + b\tau, \tau) \\
&= -\frac{1}{q} e^{-2\pi i(z+(b-1)\tau)} \theta(z + (b-1)\tau, \tau) \\
&= -\frac{1}{q} e^{-2\pi i(z+(b-1)\tau)} \left(-\frac{1}{q}\right) e^{-2\pi i(z+(b-2)\tau)} \theta(z + (b-2)\tau, \tau) \\
&= (-1)^b \frac{1}{q^b} e^{-2\pi i[z+(b-1)\tau)+(z+(b-2)\tau)+\cdots+z]} \theta(z, \tau) \\
&= (-1)^b \frac{1}{q^b} e^{-2\pi i b z - \pi i b(b-1)\tau} \theta(z, \tau) \\
&= (-1)^b e^{-2\pi i b z - \pi i b^2 \tau} \theta(z, \tau).
\end{aligned}
\tag{2.15}
$$

Similarly,

$$\theta_1(z + a, \tau) = (-1)^a \theta_1(z, \tau), \ \theta_1(z + b\tau, \tau) = e^{-2\pi i b z - \pi i b^2 \tau} \theta_1(z, \tau); \tag{2.16}$$

$$\theta_2(z + a, \tau) = \theta_2(z, \tau), \ \theta_2(z + b\tau, \tau) = (-1)^b e^{-2\pi i b z - \pi i b^2 \tau} \theta_2(z, \tau); \tag{2.17}$$

$$\theta_3(z + a, \tau) = \theta_3(z, \tau), \ \theta_3(z + b\tau, \tau) = e^{-2\pi i b z - \pi i b^2 \tau} \theta_3(z, \tau). \tag{2.18}$$

**Definition 2.1.** *Let $\Gamma$ be a subgroup of $SL(2, \mathbf{Z})$. A modular form over $\Gamma$ is a holomorphic function $f(\tau)$ on $\mathbf{H}$ such that for any*

$$g = \begin{pmatrix} a_1 & a_2 \\ a_3 & a_4 \end{pmatrix} \in \Gamma,$$

*the following property holds,*

$$f(g\tau) := f\left(\frac{a_1 \tau + a_2}{a_3 \tau + a_4}\right) = \chi(g)(a_3 \tau + a_4)^k f(\tau),$$

*where $\chi : \Gamma \to \mathbf{C}^*$ is a character of $\Gamma$ and $k$ is called the weight of $f$.*

2.2. **Chern-Weil forms of some characteristic classes.** Let $M$ be a smooth Riemannian manifold. Let $\nabla^{TM}$ be the associated Levi-Civita connection on $TM$ and $R^{TM} = (\nabla^{TM})^2$ be the curvature of $\nabla^{TM}$. Then $\nabla^{TM}$ extends canonically to a Hermitian connection $\nabla^{T_\mathbf{C} M}$ on $T_\mathbf{C} M = TM \otimes \mathbf{C}$.

Let $\widehat{A}(TM, \nabla^{TM})$ be the Hirzebruch $\widehat{A}$-form defined by (cf. [26])

$$\widehat{A}(TM, \nabla^{TM}) = \det{}^{1/2}\left(\frac{\frac{\sqrt{-1}}{4\pi} R^{TM}}{\sinh\left(\frac{\sqrt{-1}}{4\pi} R^{TM}\right)}\right). \tag{2.19}$$

Let $E$, $F$ be two Hermitian vector bundles over $M$ carrying Hermitian connections $\nabla^E$, $\nabla^F$ respectively. Let $R^E = (\nabla^E)^2$ (resp. $R^F = (\nabla^F)^2$) be the curvature of $\nabla^E$ (resp. $\nabla^F$). If we set the formal difference $G = E - F$, then $G$ carries an induced Hermitian connection $\nabla^G$ in an obvious sense. We define the associated Chern character form as (cf. [26])

$$\mathrm{ch}(G, \nabla^G) = \mathrm{tr}\left[\exp\left(\frac{\sqrt{-1}}{2\pi} R^E\right)\right] - \mathrm{tr}\left[\exp\left(\frac{\sqrt{-1}}{2\pi} R^F\right)\right]. \tag{2.20}$$

For any complex number $t$, let

$$S_t(E) = \mathbf{C}|_M + tE + t^2 S^2(E) + \cdots, \ \Lambda_t(E) = \mathbf{C}|_M + tE + t^2 \Lambda^2(E) + \cdots$$

denote respectively the total symmetric and exterior powers of $E$, which lie in $K(M)[[t]]$. The following relations between these two operations hold (cf. [2]),

$$(2.21) \qquad S_t(E) = \frac{1}{\Lambda_{-t}(E)}, \qquad \Lambda_t(E - F) = \frac{\Lambda_t(E)}{\Lambda_t(F)}.$$

The connections $\nabla^E$ and $\nabla^F$ naturally induce connections on $S_t(E)$ and $\Lambda_t(E)$ etc. Moreover, if $\{\omega_i\}$, $\{\omega_j'\}$ are formal Chern roots for Hermitian vector bundles $E$, $F$ respectively, then [13]

$$(2.22) \qquad \operatorname{ch}\left(\Lambda_t(E), \nabla^{\Lambda_t(E)}\right) = \prod_i (1 + e^{\omega_i} t).$$

Therefore, we have the following formulas for Chern character forms,

$$(2.23) \qquad \operatorname{ch}\left(S_t(E), \nabla^{S_t(E)}\right) = \frac{1}{\operatorname{ch}\left(\Lambda_{-t}(E), \nabla^{\Lambda_{-t}(E)}\right)} = \frac{1}{\prod_i (1 - e^{\omega_i} t)},$$

$$(2.24) \qquad \operatorname{ch}\left(\Lambda_t(E - F), \nabla^{\Lambda_t(E-F)}\right) = \frac{\operatorname{ch}\left(\Lambda_t(E), \nabla^{\Lambda_t(E)}\right)}{\operatorname{ch}\left(\Lambda_t(F), \nabla^{\Lambda_t(F)}\right)} = \frac{\prod_i (1 + e^{\omega_i} t)}{\prod_j (1 + e^{\omega_j'} t)}.$$

If $V$ is a real Euclidean vector bundle over $M$ carrying a Euclidean connection $\nabla^V$, then its complexification $V_{\mathbf{C}} = V \otimes \mathbf{C}$ is a complex vector bundle over $M$ carrying a canonically induced Hermitian metric from that of $V$, as well as a Hermitian connection $\nabla^{V_{\mathbf{C}}}$ induced from $\nabla^V$. If $E$ is a complex vector bundle over $M$, set $\widetilde{E} = E - \mathbf{C}^{\operatorname{rk}(E)} \in K(M)$.

Let $\Theta(T_{\mathbf{C}} M)$ be the Witten bundle over $M$ defined by $\Theta(T_{\mathbf{C}} M) := \overset{\infty}{\underset{m=1}{\otimes}} S_{q^{2m}}(\widetilde{T_{\mathbf{C}} M})$. The connection $\nabla^{TM}$ induces a connection $\nabla^{\Theta(T_{\mathbf{C}} M)}$ on $\Theta(T_{\mathbf{C}} M)$. Then the *Witten form*

$$\mathcal{W}(M, \nabla^{TM}) := \widehat{A}(TM, \nabla^{TM}) \operatorname{ch}\left(\Theta\left(T_{\mathbf{C}} M\right), \nabla^{\Theta(T_{\mathbf{C}} M)}\right)$$

can be expressed in terms of the theta function via curvature as $\det^{\frac{1}{2}}\left(\frac{R^{TM}}{4\pi^2} \frac{\theta'(0, \tau)}{\theta\left(\frac{R^{TM}}{4\pi^2}, \tau\right)}\right)$ (cf. [6]).

If $M$ is oriented, then

$$(2.25) \qquad W(M) = \int_M \mathcal{W}(M, \nabla^{TM}) = \int_M \det^{\frac{1}{2}}\left(\frac{R^{TM}}{4\pi^2} \frac{\theta'(0, \tau)}{\theta\left(\frac{R^{TM}}{4\pi^2}, \tau\right)}\right).$$

By using the Chern root algorithm, the Witten genus can be expressed as in (1.1).

## 3. Generalized Witten genera

In this section, we discuss two extensions of the Witten genus. One concerns spin$^c$ manifolds and the other one is a mod 2 extension for $8k + 2$ spin manifolds.

Let $M$ be an even dimensional closed oriented spin$^c$-manifold and $L$ be the complex line bundle associated to the given spin$^c$ structure on $M$ (cf. [17, Appendix]). Denote by $c = c_1(L)$ the first Chern class of $L$. Also, we use $L_{\mathbf{R}}$ for the notation of $L$, when it is viewed as an oriented real plane bundle.

Let $\Theta(T_{\mathbf{C}}M, L_{\mathbf{R}} \otimes \mathbf{C})$, $\Theta^*(T_{\mathbf{C}}M, L_{\mathbf{R}} \otimes \mathbf{C})$ be the virtual complex vector bundles over $M$ defined by

$$\Theta(T_{\mathbf{C}}M, L_{\mathbf{R}} \otimes \mathbf{C}) := \left( \overset{\infty}{\underset{m=1}{\otimes}} S_{q^{2m}}(\widetilde{T_{\mathbf{C}}M}) \right) \otimes \left( \overset{\infty}{\underset{n=1}{\otimes}} \Lambda_{q^{2n}}(\widetilde{L_{\mathbf{R}} \otimes \mathbf{C}}) \right)$$
$$\otimes \left( \overset{\infty}{\underset{u=1}{\otimes}} \Lambda_{-q^{2u-1}}(\widetilde{L_{\mathbf{R}} \otimes \mathbf{C}}) \right) \otimes \left( \overset{\infty}{\underset{v=1}{\otimes}} \Lambda_{q^{2v-1}}(\widetilde{L_{\mathbf{R}} \otimes \mathbf{C}}) \right),$$

and

$$\Theta^*(T_{\mathbf{C}}M, L_{\mathbf{R}} \otimes \mathbf{C}) := \left( \overset{\infty}{\underset{m=1}{\otimes}} S_{q^{2m}}(\widetilde{T_{\mathbf{C}}M}) \right) \otimes \left( \overset{\infty}{\underset{n=1}{\otimes}} \Lambda_{-q^{2n}}(\widetilde{L_{\mathbf{R}} \otimes \mathbf{C}}) \right).$$

Let $g^{TM}$ be a Riemannian metric on $M$. Let $\nabla^{TM}$ be the Levi-Civita connection associated to $g^{TM}$. Let $g^{T_{\mathbf{C}}M}$ and $\nabla^{T_{\mathbf{C}}M}$ be the induced Hermitian metric and Hermitian connection on $T_{\mathbf{C}}M$. Let $h^L$ be a Hermitian metric on $L$ and $\nabla^L$ be a Hermitian connection. Let $h^{L_{\mathbf{R}}}$ and $\nabla^{L_{\mathbf{R}}}$ be the induced Euclidean metric and connection on $L_{\mathbf{R}}$. Then $\nabla^{TM}$ and $\nabla^L$ induce connections $\nabla^{\Theta(T_{\mathbf{C}}M, L_{\mathbf{R}} \otimes \mathbf{C})}$ and $\nabla^{\Theta^*(T_{\mathbf{C}}M, L_{\mathbf{R}} \otimes \mathbf{C})}$ on $\Theta(T_{\mathbf{C}}M, L_{\mathbf{R}} \otimes \mathbf{C})$ and $\Theta^*(T_{\mathbf{C}}M, L_{\mathbf{R}} \otimes \mathbf{C})$ respectively. Let $c$ be the first Chern form of $(L, \nabla^L)$.

Define the *generalized Witten forms*

$$\mathcal{W}_c(M) := \widehat{A}(TM, \nabla^{TM}) \exp\left(\frac{c}{2}\right) \operatorname{ch}\left(\Theta(T_{\mathbf{C}}M, L_{\mathbf{R}} \otimes \mathbf{C}), \nabla^{\Theta(T_{\mathbf{C}}M, L_{\mathbf{R}} \otimes \mathbf{C})}\right)$$

if $\dim M = 4k$ and

$$\mathcal{W}_c(M) := \widehat{A}(TM, \nabla^{TM}) \exp\left(\frac{c}{2}\right) \operatorname{ch}\left(\Theta^*(T_{\mathbf{C}}M, L_{\mathbf{R}} \otimes \mathbf{C}), \nabla^{\Theta^*(T_{\mathbf{C}}M, L_{\mathbf{R}} \otimes \mathbf{C})}\right)$$

if $\dim M = 4k + 2$.

One can express the generalized Witten forms in terms of theta functions via curvatures by

$$\mathcal{W}_c(M^{4k}) = \det^{\frac{1}{2}}\left( \frac{R^{TM}}{4\pi^2} \frac{\theta'(0, \tau)}{\theta\left(\frac{R^{TM}}{4\pi^2}, \tau\right)} \right) \frac{\theta_1\left(\frac{R^L}{4\pi^2}, \tau\right)}{\theta_1(0, \tau)} \frac{\theta_2\left(\frac{R^L}{4\pi^2}, \tau\right)}{\theta_2(0, \tau)} \frac{\theta_3\left(\frac{R^L}{4\pi^2}, \tau\right)}{\theta_3(0, \tau)}$$

and

$$\mathcal{W}_c(M^{4k+2}) = \det^{\frac{1}{2}}\left( \frac{R^{TM}}{4\pi^2} \frac{\theta'(0, \tau)}{\theta\left(\frac{R^{TM}}{4\pi^2}, \tau\right)} \right) \frac{\sqrt{-1}\theta\left(\frac{R^L}{4\pi^2}, \tau\right)}{\theta_1(0, \tau)\theta_2(0, \tau)\theta_3(0, \tau)},$$

where $R^{TM} = (\nabla^{TM})^2$ and $R^L = (\nabla^L)^2$ are the curvatures.

We can also express the generalized Witten forms by using the Chern-root algorithm. Let $\{\pm 2\pi\sqrt{-1}z_j\}$ be the formal Chern roots for $(T_{\mathbf{C}}M, \nabla^{T_{\mathbf{C}}M})$ and set $u = \frac{1}{2\pi\sqrt{-1}}c$. In terms of the theta-functions, we get through direct computations that

$$(3.1) \qquad \mathcal{W}_c(M^{4k}) = \left( \prod_{j=1}^{2k} z_j \frac{\theta'(0, \tau)}{\theta(z_j, \tau)} \right) \frac{\theta_1(u, \tau)\theta_2(u, \tau)\theta_3(u, \tau)}{\theta_1(0, \tau)\theta_2(0, \tau)\theta_3(0, \tau)},$$

and

$$(3.2) \qquad \mathcal{W}_c(M^{4k+2}) = \left( \prod_{j=1}^{2k+1} z_j \frac{\theta'(0, \tau)}{\theta(z_j, \tau)} \right) \frac{\sqrt{-1}\theta(u, \tau)}{\theta_1(0, \tau)\theta_2(0, \tau)\theta_3(0, \tau)}.$$

Define the *generalized Witten genus* by

$$(3.3) \qquad W_c(M^{4k}) := \int_{M^{4k}} \mathcal{W}_c(M^{4k}),$$

and
$$W_c(M^{4k+2}) := \int_{M^{4k+2}} \mathcal{W}_c(M^{4k+2}). \tag{3.4}$$

Let $S_c(TM) = S_{c,+}(TM) \oplus S_{c,-}(TM)$ denote the bundle of spinors associated to the spin$^c$ structure, $(TM, g^{TM})$ and $(L, h^L)$. Then $S_c(TM)$ carries induced Hermitian metric and connection preserving the above $\mathbf{Z}_2$-grading. Let $D_{c,\pm} : \Gamma(S_{c,\pm}(TM)) \to \Gamma(S_{c,\mp}(TM))$ denote the induced spin$^c$ Dirac operators (cf. [17]). Consider the (virtual) complex vector bundles

$$\Theta(T_{\mathbf{C}}M, L_{\mathbf{R}} \otimes \mathbf{C}) \oplus L^{-1} \otimes \Theta(T_{\mathbf{C}}M, L_{\mathbf{R}} \otimes \mathbf{C}), \tag{3.5}$$

$$\Theta^*(T_{\mathbf{C}}M, L_{\mathbf{R}} \otimes \mathbf{C}) \ominus L^{-1} \otimes \Theta^*(T_{\mathbf{C}}M, L_{\mathbf{R}} \otimes \mathbf{C}). \tag{3.6}$$

They carry induced Hermitian metrics and connections. By the Atiyah-Singer index theorem for spin$^c$ manifolds, it's not hard to see that

$$W_c(M^{4k}) = \frac{1}{2}\mathrm{Ind}\left(D_{c,+} \otimes \left(\Theta(T_{\mathbf{C}}M, L_{\mathbf{R}} \otimes \mathbf{C}) \oplus L^{-1} \otimes \Theta(T_{\mathbf{C}}M, L_{\mathbf{R}} \otimes \mathbf{C})\right)\right) \in \mathbf{Z}[[q]]$$

and

$$W_c(M^{4k+2}) = \frac{1}{2}\mathrm{Ind}\left(D_{c,+} \otimes \left(\Theta^*(T_{\mathbf{C}}M, L_{\mathbf{R}} \otimes \mathbf{C}) \ominus L^{-1} \otimes \Theta^*(T_{\mathbf{C}}M, L_{\mathbf{R}} \otimes \mathbf{C})\right)\right) \in \mathbf{Z}[[q]].$$

We call the operators

$$D_{c,+} \otimes \left(\Theta(T_{\mathbf{C}}M, L_{\mathbf{R}} \otimes \mathbf{C}) \oplus L^{-1} \otimes \Theta(T_{\mathbf{C}}M, L_{\mathbf{R}} \otimes \mathbf{C})\right)$$

and

$$D_{c,+} \otimes \left(\Theta^*(T_{\mathbf{C}}M, L_{\mathbf{R}} \otimes \mathbf{C}) \ominus L^{-1} \otimes \Theta^*(T_{\mathbf{C}}M, L_{\mathbf{R}} \otimes \mathbf{C})\right)$$

*the generalized Witten operators.*

Clearly, when the spin$^c$ manifold $M$ is actually spin with $c = 0$, $W_c(M)$ is exactly the Witten genus of $M$.

Since $M$ is spin$^c$, $TM \oplus L_{\mathbf{R}}$ is spin. By a result of McLaughlin ([20], Lemma 2.2), there is a class $\lambda_c \in H^4(M, \mathbf{Z})$ associated to the *spin$^c$* structure such that $2\lambda_c = p_1(TM \oplus L_{\mathbf{R}})$. However,

$$\begin{aligned}
p_1(TM \oplus L_{\mathbf{R}}) \\
&= -c_2((TM \oplus L_{\mathbf{R}}) \otimes \mathbf{C}) \\
&= -c_2(T_{\mathbf{C}}M) - c_2(L \oplus \overline{L}) - c_1(T_{\mathbf{C}}M)c_1(L \oplus \overline{L}) \\
&= p_1(TM) + c^2.
\end{aligned} \tag{3.7}$$

So $2\lambda_c = p_1(TM) + c^2$.

The following result can be proved by a direct computation.

**Theorem 3.1.** *(i) If $\dim M = 4k$ and $\lambda_c - 2c^2 = 0$ (it is enough to assume that $p_1(TM) - 3c^2 = 0$ rationally), then $W_c(M)$ is an integral modular form of weight $2k$ over $SL(2, \mathbf{Z})$;*

*(ii) If $\dim M = 4k+2$ and $\lambda_c - c^2 = 0$ (it is enough to assume that $p_1(TM) - c^2 = 0$ rationally), then $W_c(M)$ is an integral modular form of weight $2k$ over $SL(2, \mathbf{Z})$.*

Simply denote the class $\lambda_c - 2c^2$ by $\frac{p_1(TM) - 3c^2}{2}$ and $\lambda_c - c^2$ by $\frac{p_1(TM) - c^2}{2}$. Inspired by Theorem 3.1, we made in [8] the following definition:

**Definition 3.1.** *If $M$ is an even dimensional closed oriented spin$^c$ manifold and $L$ denotes the complex line bundle associated to a given spin$^c$ structure on $M$ with $c = c_1(L)$ denotes the first Chern class of $L$. Then $M$ is called a string$^c$ manifold if $\frac{p_1(TM) - (2 + (-1)^{\frac{\dim M}{2}})c^2}{2} = 0$.*

With this terminology, Theorem 3.1 can then be stated simply as follows: the genus $W_c(M)$ defined in (3.3) and (3.4) is an integral modular form on a string$^c$ manifold.

As a special example, when $M$ is 4 dimensional,

$$(3.8) \qquad W(M) = \int_M \frac{z_1 \theta'(0,\tau)}{\theta(z_1,\tau)} \frac{z_2 \theta'(0,\tau)}{\theta(z_2,\tau)} \frac{\theta_1(u,\tau)\theta_2(u,\tau)\theta_3(u,\tau)}{\theta_1(0,\tau)\theta_2(0,\tau)\theta_3(0,\tau)}.$$

Up to the terms of degree higher than 2, we have

$$(3.9) \qquad \frac{z_j}{\theta(z_j,\tau)} = \frac{1}{2\pi q^{1/4} \prod_{l=1}^{\infty}[(1-q^{2l})^3]} \left[1 + \left(\frac{1}{6}\pi^2 - \sum_{i=1}^{\infty} \frac{4\pi^2 q^{2i}}{(1-q^{2i})^2}\right) z_j^2\right], \ j=1,2$$

and

$$(3.10) \qquad \theta_1(u,\tau)\theta_2(u,\tau)\theta_3(u,\tau) = 2q^{1/4} \prod_{l=1}^{\infty}[(1-q^{2l})^3] \left[1 + \left(-\frac{1}{2}\pi^2 + \sum_{i=1}^{\infty} \frac{12\pi^2 q^{2i}}{(1-q^{2i})^2}\right) u^2\right].$$

Thus up to some constant, we have

$$(3.11)$$
$$W(M) = \int_M \left(\frac{1}{6}\pi^2 - \sum_{i=1}^{\infty} \frac{4\pi^2 q^{2i}}{(1-q^{2i})^2}\right)(z_1^2 + z_2^2 - 3u^2) = \left(\frac{1}{6}\pi^2 - \sum_{i=1}^{\infty} \frac{4\pi^2 q^{2i}}{(1-q^{2i})^2}\right) \frac{1}{-4\pi^2} (p_1 - 3c^2).$$

Therefore one can see that *if $M$ is a 4-dimensional string$^c$ manifold, then $W_c(M) = 0$.*

The homotopy theoretical refinements of the classical Witten genus leads to the theory of topological modular forms developed by Hopkins and Miller [15]. In [1], it is shown that the Witten genus

$$\pi_* MString \to MF_*$$

is the value on homotopy group of a map of $E_\infty$-spectra

$$MString \to \text{tmf}.$$

It would be interesting to study the homotopy theoretical refinements of the generalized Witten genus for string$^c$ manifolds.

Similar to the Witten operator, the generalized Witten operators are also rigid. For any integer $n \geq 0$, let $R_n(T_\mathbf{C}M, L_\mathbf{R} \otimes \mathbf{C})$, $Q_n(T_\mathbf{C}M, L_\mathbf{R} \otimes \mathbf{C})$ be the (virtual) complex vector bundles defined by

$$(3.12) \qquad \Theta(T_\mathbf{C}M, L_\mathbf{R} \otimes \mathbf{C}) \oplus L^{-1} \otimes \Theta(T_\mathbf{C}M, L_\mathbf{R} \otimes \mathbf{C}) = \sum_{n=0}^{+\infty} R_n(T_\mathbf{C}M, L_\mathbf{R} \otimes \mathbf{C}) q^n,$$

$$(3.13) \qquad \Theta^*(T_\mathbf{C}M, L_\mathbf{R} \otimes \mathbf{C}) \ominus L^{-1} \otimes \Theta^*(T_\mathbf{C}M, L_\mathbf{R} \otimes \mathbf{C}) = \sum_{n=0}^{+\infty} Q_n(T_\mathbf{C}M, L_\mathbf{R} \otimes \mathbf{C}) q^{2n}$$

respectively. Clearly, each of the $R_n(T_\mathbf{C}M, L_\mathbf{R} \otimes \mathbf{C})$'s and $Q_n(T_\mathbf{C}M, L_\mathbf{R} \otimes \mathbf{C})$'s carries induced Hermitian metric and connection.

Now we assume that $M$ admits an $S^1$-action which lifts to an action on $L$. Moreover, we assume that this action preserves the given spin$^c$-structure associated to $(M,L)$, as well as the metrics and connections involved.

For any integer $n \geq 0$, let $D_{c,\pm}^{R_n(T_\mathbf{C}M, L_\mathbf{R} \otimes \mathbf{C})} : \Gamma(S_{c,\pm}(TM) \otimes R_n(T_\mathbf{C}M, L_\mathbf{R} \otimes \mathbf{C})) \to \Gamma(S_{c,\mp}(TM) \otimes R_n(T_\mathbf{C}M, L_\mathbf{R} \otimes \mathbf{C}))$ denote the corresponding twisted spin$^c$ Dirac operators. Let $D_{c,\pm}^{Q_n(T_\mathbf{C}M, L_\mathbf{R} \otimes \mathbf{C})}$ be defined similarly. Then each of the $D_{c,\pm}^{R_n(T_{\mathbf{C}}M, L_{\mathbf{R}} \otimes \mathbf{C})}$'s and $D_{c,\pm}^{Q_n(T_{\mathbf{C}}M, L_{\mathbf{R}} \otimes \mathbf{C})}$'s is an $S^1$-equivariant operator.

Let $ES^1$ be the universal $S^1$ principal bundle over the classifying space $BS^1$ of $S^1$. Let $H^*(BS^1, \mathbf{Z}) = \mathbf{Z}[[u]]$, with $u$ a generator of degree 2. Let $\pi : M \times_{S^1} ES^1 \to BS^1$ be the projection.

**Theorem 3.2.** (i) *Let $M$ be a spin$^c$ manifold of dimension $4k$ whose first equivariant Pontryajin class verifies $3c_1(L)_{S^1}^2 - p_1(TM)_{S^1} = l\pi^*u^2$, where $l$ is an integer and $c_1(L)_{S^1}$ is the first equivariant Chern class of $L$. One has*
a) *if $l = 0$, then for any integer $n \geq 0$, $D_{c,+}^{R_n(T_{\mathbf{C}}M, L_{\mathbf{R}} \otimes \mathbf{C})}$ is rigid in the sense of Witten [22, Section 4];*
b) *if $l < 0$, then for any integer $n \geq 0$, $\mathrm{ind}(D_{c,+}^{R_n(T_{\mathbf{C}}M, L_{\mathbf{R}} \otimes \mathbf{C})}) = 0$.*

(ii) *Let $M$ be a spin$^c$ manifold of dimension $4k+2$ whose first equivariant Pontryajin class verifies $c_1(L)_{S^1}^2 - p_1(TM)_{S^1} = l\pi^*u^2$, where $l$ is an integer and $c_1(L)_{S^1}$ is the first equivariant Chern class of $L$. One has*
a) *if $l = 0$, then for any integer $n \geq 0$, $D_{c,+}^{Q_n(T_{\mathbf{C}}M, L_{\mathbf{R}} \otimes \mathbf{C})}$ is rigid in the sense of Witten [22, Section 4];*
b) *if $l < 0$, then for any integer $n \geq 0$, $\mathrm{ind}(D_{c,+}^{Q_n(T_{\mathbf{C}}M, L_{\mathbf{R}} \otimes \mathbf{C})}) = 0$.*

*Proof.* If $c \in H^2(M, \mathbf{Z})$ is even, then $M$ is spin and the corresponding twisted Dirac operators have been constructed in [18, Page 356, Example f], [18, Page 353, Example b] and [22, (34)] respectively, and one can apply [18, Theorem 1] to obtain the required rigidity.

In the general case where $c$ needs not be even, we find that the arguments in [18] can be applied here to obtain the required rigidity. Q.E.D.

**Remark 3.1.** In view of the modularity and the rigidity associated to the formal Dirac operators $\sum_{n=0}^{+\infty} D_{c,+}^{R_n(T_{\mathbf{C}}M, L_{\mathbf{R}} \otimes \mathbf{C})} q^n$ and $\sum_{n=0}^{+\infty} D_{c,+}^{Q_n(T_{\mathbf{C}}M, L_{\mathbf{R}} \otimes \mathbf{C})} q^{2n}$ on string$^c$ manifolds, one would expect that, just like in the string manifolds case considered in [22], in general these operators should be viewed as a kind of Dirac operators on loop spaces of spin$^c$ manifolds.

Another extension of the classical Witten genus concerns $8k+2$-dimensional spin manifolds. Let $B$ be an $8k+2$ dimensional compact spin manifold. Define
$$\phi(B) := \mathrm{ind}_2(\Theta(TB)) \in \mathbf{Z}_2[[q]],$$
where $\mathrm{ind}_2$ is the mod 2 index in the sense of Atiyah and Singer [4] and $\Theta(TB) = \overset{\infty}{\underset{m=1}{\otimes}} S_{q^{2m}}(\widetilde{TB})$ is the Witten bundle constructed in [24, (17)][3].

It's not hard to see that $\phi$ defines a ring homomorphism
$$(3.14) \qquad \phi : \Omega_{8k+2}^{\mathrm{spin}} \to KO^{-(8k+2)}(\mathrm{pt.})[[q]].$$

We call it the *mod 2 Witten genus*. It would be very interesting to know if $\phi(B)$ is a level 1 modular form over the ring $KO^{-(8k+2)}(\mathrm{pt.})$ when $B$ is string [4]. We will study this question in further article.

**Remark 3.2.** The construction of $W_c(M)$ is actually inspired by the mod 2 Witten genus $\phi(B)$ considered here, as well as the considerations in [12] where the mod 2 elliptic genera are studied in the framework of Rokhlin congruences.

---

[3]Or rather its real form.

[4]In some sense, $\phi$ might be thought of as an analogue of Ochanine's $\beta$ invariant defined in [21], which is a level 2 modular form over the ring $KO^{-(8k+2)}(\mathrm{pt.})$ on any compact spin manifold.

## 4. Vanishing theorems on generalized complete intersections

In this section, we give the Landweber-Stong type vanishing theorems on the generalized Witten genus and the mod 2 Witten genus constructed in Section 3.

### 4.1. Vanishing theorems.

We first recall the result of Chen-Han in [7] for completeness. Let $V_{(d_{\alpha\beta})}$ be a nonsingular $4k$ dimensional generalized complete intersection in the product of complex projective spaces $\mathbf{C}P^{n_1} \times \mathbf{C}P^{n_2} \times \cdots \times \mathbf{C}P^{n_s}$, such that $[V_{(d_{\alpha\beta})}] \in H_{4k}(\mathbf{C}P^{n_1} \times \mathbf{C}P^{n_2} \times \cdots \times \mathbf{C}P^{n_s}, \mathbf{Z})$ is dual to $\prod_{\alpha=1}^{t}(\sum_{\beta=1}^{s} d_{\alpha\beta}x_\beta) \in H^{2t}(\mathbf{C}P^{n_1} \times \mathbf{C}P^{n_2} \times \cdots \times \mathbf{C}P^{n_s}, \mathbf{Z})$, where each $x_\beta \in H^2(\mathbf{C}P^{n_\beta}, \mathbf{Z})$, $1 \leq \beta \leq s$, is the generator of $H^*(\mathbf{C}P^{n_\beta}, \mathbf{Z})$ and $d_{\alpha\beta}$, $1 \leq \alpha \leq t$, $1 \leq \beta \leq s$, are integers.

Let $P_\beta : \mathbf{C}P^{n_1} \times \mathbf{C}P^{n_2} \times \cdots \times \mathbf{C}P^{n_s} \to \mathbf{C}P^{n_\beta}$, $1 \leq \beta \leq s$, be the projection on the $\beta$-th factor. Then $V_{(d_{\alpha\beta})}$ is the transversal intersection of zero loci of smooth global sections of line bundles $\overset{s}{\underset{\beta=1}{\otimes}} P_\beta^*(\mathcal{O}_\beta(d_{\alpha\beta}))$, $1 \leq \alpha \leq t$, where $\mathcal{O}_\beta(d_{\alpha\beta})$ denotes the $d_{\alpha\beta}$-th power of the canonical line bundle $\mathcal{O}(1)$ over $\mathbf{C}P^{n_\beta}$.[5]

Putting some relevant conditions on the data $n_\beta$, $1 \leq \beta \leq s$, and $d_{\alpha\beta}$, $1 \leq \alpha \leq t$, $1 \leq \beta \leq s$, the generalized complete intersection $V_{(d_{\alpha\beta})}$ can be made string. This systematically provides us with a lot of interesting examples of string manifolds.

Set
$$D = \begin{bmatrix} d_{11} & d_{12} & \cdots & d_{1s} \\ d_{21} & d_{22} & \cdots & d_{2s} \\ \cdots & \cdots & \cdots & \cdots \\ d_{t1} & d_{t2} & \cdots & d_{ts} \end{bmatrix}.$$

Let $m_\beta$ be the number of nonzero elements in the $\beta$-th column of $D$.

The main result in [7] says *if $m_\beta + 2 \leq n_\beta$, $1 \leq \beta \leq s$, and $V_{(d_{\alpha\beta})}$ is string, then the Witten genus $W(V_{(d_{\alpha\beta})})$ vanishes.*

Putting $s = 1$, one obtains the Landweber-Stong vanishing theorem: *the Witten genus vanishes on string complete intersections in a complex projective space.*

It is natural to ask if there are Landweber-Stong type vanishing theorems for the generalized Witten genus and the mod 2 Witten genus constructed in Section 3. We provide the following answers to this natural question.

**Theorem 4.1.** *If $m_\beta + 2 \leq n_\beta$, $1 \leq \beta \leq s$, and $V_{(d_{\alpha\beta})}$ is string$^c$ with $L = \mathfrak{L}|_{V_{(d_{\alpha\beta})}}$ for some line bundle $\mathfrak{L}$ on $\mathbf{C}P^{n_1} \times \mathbf{C}P^{n_2} \times \cdots \times \mathbf{C}P^{n_s}$, then the genus $W_c(V_{(d_{\alpha\beta})})$ vanishes.*

Specializing to the case $s = 1$, i.e., we only consider the complete intersection $V(d_1, d_2, \cdots, d_r)$ in $\mathbf{C}P^{k+r}$ (the conditions that $d_i > 0$, $0 \leq i \leq r$, are also put). In this case, the conditions in Theorem 4.1 can be simplified thanks to the well known Lefschetz hyperplane theorem on nonsingular complex projective algebraic varieties (cf. Chapter 1 in [10]). In fact, the Lefschetz hyperplane theorem asserts that
$$i^* : H^j(\mathbf{C}P^{k+r}, \mathbf{Z}) \to H^j(V(d_1, d_2, \cdots, d_r), \mathbf{Z})$$
is an isomorphism when $j < k$. Therefore from Theorem 4.1, one sees that (taking $j = 2$) the genus $W_c$ vanishes on string$^c$ complete intersections of complex dimension $\geq 3$ in a complex projective space.

---

[5]In [9], Gorbounov and Ochanine calculated the complex elliptic genus for such generalized complete intersections and showed that it equals to the elliptic genus of the Landau-Ginzburg model, which are the mirror partners of these generalized complete intersections according to Hori and Vafa [16].

On the other hand, we have already shown that $W_c$ vanishes on string$^c$ manifolds of real dimension 4 in Section 3. Therefore we obtain the following generalization of the vanishing result of Landweber-Stong mentioned before.

**Corollary 4.1.** *The genus $W_c$ vanishes on string$^c$ complete intersections of complex dimension $k$ ($k \geq 2$) in a complex projective space.*

In the mod 2 case, we have the following vanishing theorem:

**Theorem 4.2.** *If $m_\beta + 2 \leq n_\beta, 1 \leq \beta \leq s$, $\dim_{\mathbf{R}} V_{(d_{\alpha\beta})} = 8k + 2$, $V_{(d_{\alpha\beta})}$ is string and one of the generalized hypersurfaces that generate $V_{(d_{\alpha\beta})}$ is of even degree (i.e. one of the rows of D consists of only even numbers), then the mod 2 Witten genus $\phi(V_{(d_{\alpha\beta})})$ vanishes.*

Putting $s = 1$, we obtain the following corollary (see the deduction of Corollary 4.2 from Theorem 4.2 at the end of Section 4.5), which is the mod 2 analogue of the Landweber-Stong vanishing result.

**Corollary 4.2.** *If $V(d_1, \cdots, d_r) \subset \mathbf{C}P^{4k+1+r}$ is a string complete intersection of complex dimension $4k+1$ ($k \geq 1$) in a complex projective space, then the mod 2 Witten genus $\phi(V(d_1, \cdots, d_r))$ vanishes.*

**Remark 4.1.** The above Landweber-Stong type vanishing theorems and their corollaries indicate that the genera we introduce in Section 3 are in some sense the "right" extensions of the classical Witten genus.

### 4.2. Numerical characterizations of string$^c$ complete intersections. Let

$$i : V_{(d_{\alpha\beta})} \to \mathbf{C}P^{n_1} \times \mathbf{C}P^{n_2} \times \cdots \times \mathbf{C}P^{n_s}$$

denote the inclusion embedding. It's not hard to see that

$$(4.1) \qquad i^* T_{\mathbf{R}}(\mathbf{C}P^{n_1} \times \mathbf{C}P^{n_2} \times \cdots \times \mathbf{C}P^{n_s}) \cong TV_{(d_{\alpha\beta})} \oplus i^* \left( \bigoplus_{\alpha=1}^{t} \left( \bigotimes_{\beta=1}^{s} P_\beta^*(\mathcal{O}_\beta(d_{\alpha\beta})) \right) \right),$$

where we forget the complex structure of the line bundles $\bigotimes_{\beta=1}^{s} P_\beta^*(\mathcal{O}_\beta(d_{\alpha\beta}))$, $1 \leq \alpha \leq t$. Therefore for the total Stiefel-Whitney class, we have

$$(4.2) \qquad i^* w(T_{\mathbf{R}}(\mathbf{C}P^{n_1} \times \mathbf{C}P^{n_2} \times \cdots \times \mathbf{C}P^{n_s})) = w(TV_{(d_{\alpha\beta})}) \prod_{\alpha=1}^{t} i^* w(\bigotimes_{\beta=1}^{s} P_\beta^*(\mathcal{O}_\beta(d_{\alpha\beta}))),$$

or more precisely,

$$(4.3) \qquad i^* \left( \prod_{\beta=1}^{s} (1+x_\beta)^{n_\beta+1} \right) \equiv w(TV_{(d_{\alpha\beta})}) \prod_{\alpha=1}^{t} i^* \left( 1 + \sum_{\beta=1}^{s} d_{\alpha\beta} x_\beta \right) \mod 2,$$

where each $x_\beta \in H^2(\mathbf{C}P^{n_\beta}, \mathbf{Z})$ is the generator of $H^*(\mathbf{C}P^{n_\beta}, \mathbf{Z})$.

By (4.3), we see that

$$(4.4) \qquad\qquad\qquad w_1(TV_{(d_{\alpha\beta})}) = 0,$$

$$(4.5) \qquad\qquad w_2(TV_{(d_{\alpha\beta})}) \equiv \sum_{\beta=1}^{s} \left( n_\beta + 1 - \sum_{\alpha=1}^{t} d_{\alpha\beta} \right) i^* x_\beta \mod 2.$$

As for the total rational Pontryagin class, we have

$$(4.6) \qquad i^* p(T_{\mathbf{R}}(\mathbf{C}P^{n_1} \times \mathbf{C}P^{n_2} \times \cdots \times \mathbf{C}P^{n_s})) = p(TV_{(d_{\alpha\beta})}) \prod_{\alpha=1}^{t} i^* p(\bigotimes_{\beta=1}^{s} P_\beta^*(\mathcal{O}_\beta(d_{\alpha\beta}))),$$

or

$$p(V_{(d_{\alpha\beta})}) = \prod_{\beta=1}^{s}(1+(i^*x_\beta)^2)^{n_\beta+1} \prod_{\alpha=1}^{t}\left(1+\left(\sum_{\beta=1}^{s}d_{\alpha\beta}i^*x_\beta\right)^2\right)^{-1}. \quad (4.7)$$

Hence we have

$$\begin{aligned} p_1(V_{(d_{\alpha\beta})}) &= \sum_{\beta=1}^{s}(n_\beta+1)(i^*x_\beta)^2 - \sum_{\alpha=1}^{t}\left(\sum_{\beta=1}^{s}d_{\alpha\beta}i^*x_\beta\right)^2 \\ &= \sum_{\beta=1}^{s}(n_\beta+1-\sum_{\alpha=1}^{t}d_{\alpha\beta}^2)(i^*x_\beta)^2 - \sum_{1\leq\gamma,\delta\leq s,\gamma\neq\delta}\left(\sum_{p=1}^{t}d_{\alpha\gamma}d_{\alpha\delta}(i^*x_\gamma)(i^*x_\delta)\right). \end{aligned} \quad (4.8)$$

Let $c \in H^2(\mathbf{C}P^{n_1}\times\mathbf{C}P^{n_2}\times\cdots\times\mathbf{C}P^{n_s},\mathbf{Z})$ such that $c = \sum_{\beta=1}^{s}c_\beta x_\beta$ with $c_\beta \in \mathbf{Z}$.

If $p_1(V_{(d_{\alpha\beta})}) = 3(i^*c)^2$, then

$$\begin{aligned} 0 &= i_!(p_1(V_{(d_{\alpha\beta})}) - 3(i^*c)^2) \\ &= i_!i^*\left(\sum_{\beta=1}^{s}(n_\beta+1-\sum_{\alpha=1}^{t}d_{\alpha\beta}^2)x_\beta^2 - \sum_{1\leq\gamma,\delta\leq s,\gamma\neq\delta}\left(\sum_{\alpha=1}^{t}d_{\alpha\gamma}d_{\alpha\delta}x_\gamma x_\delta\right) - 3\left(\sum_{\beta=1}^{s}c_\beta x_\beta\right)^2\right), \end{aligned}$$

i.e.,

$$\left(\prod_{\alpha=1}^{t}\left(\sum_{\beta=1}^{s}d_{\alpha\beta}x_\beta\right)\right)$$
$$\cdot\left(\sum_{\beta=1}^{s}(n_\beta+1-\sum_{\alpha=1}^{t}d_{\alpha\beta}^2)x_\beta^2 - \sum_{1\leq\gamma,\delta\leq s,\gamma\neq\delta}\left(\sum_{\alpha=1}^{t}d_{\alpha\gamma}d_{\alpha\delta}x_\gamma x_\delta\right) - 3\left(\sum_{\beta=1}^{s}c_\beta x_\beta\right)^2\right) = 0$$

in $H^{2t+4}(\mathbf{C}P^{n_1}\times\mathbf{C}P^{n_2}\times\cdots\times\mathbf{C}P^{n_s},\mathbf{Z})$, where

$$i_! : H^*(V_{(d_{\alpha\beta})},\mathbf{Z}) \to H^{*+2t}(\mathbf{C}P^{n_1}\times\mathbf{C}P^{n_2}\times\cdots\times\mathbf{C}P^{n_s},\mathbf{Z})$$

is the cohomology push forward.

Recall that $m_\beta$ is the number of nonzero elements in the $\beta$-th column of the matrix $D$ defined in Section 1, i.e., it is the number of nonzero elements in $\{d_{1\beta},\ldots,d_{t\beta}\}$.

Recall also that for any $1\leq\beta\leq s$, $1$, $x_\beta$, $x_\beta^2,\ldots,x_\beta^{n_\beta}$ are linearly independent in $H^*(\mathbf{C}P^{n_\beta},\mathbf{Z})$.

If $m_\beta + 2 \leq n_\beta$ holds for any $1\leq\beta\leq s$, then the left hand side of the above equality should be a zero polynomial of the $x_\beta$'s with $1\leq\beta\leq s$. Moreover, at least one of its factors should be zero. But $\prod_{\alpha=1}^{t}\left(\sum_{\beta=1}^{s}d_{\alpha\beta}x_\beta\right)$ is nonzero, which implies

$$\sum_{\beta=1}^{s}(n_\beta+1-\sum_{\alpha=1}^{t}d_{\alpha\beta}^2)x_\beta^2 - \sum_{1\leq\gamma,\delta\leq s,\gamma\neq\delta}(\sum_{\alpha=1}^{t}d_{\alpha\gamma}d_{\alpha\delta}x_\gamma x_\delta) = 3\left(\sum_{\beta=1}^{s}c_\beta x_\beta\right)^2$$

and consequently the following identities hold,

$$n_\beta + 1 - \sum_{\alpha=1}^{t}d_{\alpha\beta}^2 = 3c_\beta^2, \ 1\leq\beta\leq s; \quad (4.9)$$

$$(4.10) \quad -\sum_{\alpha=1}^{t} d_{\alpha\gamma}d_{\alpha\delta} = 3c_\gamma c_\delta, \ 1 \leq \gamma, \delta \leq s, \ \gamma \neq \delta.$$

In a summary, we have the following propositions.

**Proposition 4.1.** *Let $V_{(d_{\alpha\beta})}$ be a $4k$ dimensional generalized complete intersection. Let $c \in H^2(\mathbf{C}P^{n_1} \times \mathbf{C}P^{n_2} \times \cdots \times \mathbf{C}P^{n_s}, \mathbf{Z})$ such that $c = \sum_{\beta=1}^{s} c_\beta x_\beta$ with $c_\beta \in \mathbf{Z}$ and $x_\beta \in H^2(\mathbf{C}P^{n_\beta}, \mathbf{Z})$ being the generator. If $V_{(d_{\alpha\beta})}$ is string$^c$ with $c_1(L) = i^*c$ and $m_\beta + 2 \leq n_\beta$ for any $1 \leq \beta \leq s$, then the following identities hold,*

$$(4.11) \quad \begin{cases} n_\beta + 1 - \sum_{\alpha=1}^{t} d_{\alpha\beta}^2 = 3c_\beta{}^2, & 1 \leq \beta \leq s, \\ -\sum_{\alpha=1}^{t} d_{\alpha\gamma}d_{\alpha\delta} = 3c_\gamma c_\delta, & 1 \leq \gamma, \ \delta \leq s, \ \gamma \neq \delta; \end{cases}$$

*or equivalently*

$$(4.12) \quad D^T D = \mathrm{Diag}(n_1 + 1, \cdots, n_s + 1) - 3C^T C,$$

*where $C = \begin{bmatrix} c_1 & c_2 & \cdots & c_s \end{bmatrix}$.*

**Proposition 4.2.** *Let $V_{(d_{\alpha\beta})}$ be a $4k+2$ dimensional generalized complete intersection. Let $c \in H^2(\mathbf{C}P^{n_1} \times \mathbf{C}P^{n_2} \times \cdots \times \mathbf{C}P^{n_s}, \mathbf{Z})$ such that $c = \sum_{\beta=1}^{s} c_\beta x_\beta$ with $c_\beta \in \mathbf{Z}$ and $x_\beta \in H^2(\mathbf{C}P^{n_\beta}, \mathbf{Z})$ being the generator. If $V_{(d_{\alpha\beta})}$ is string$^c$ with $c_1(L) = i^*c$ and $m_\beta + 2 \leq n_\beta$ for any $1 \leq \beta \leq s$, then the following identities hold,*

$$(4.13) \quad \begin{cases} n_\beta + 1 - \sum_{\alpha=1}^{t} d_{\alpha\beta}^2 = c_\beta{}^2, & 1 \leq \beta \leq s, \\ -\sum_{\alpha=1}^{t} d_{\alpha\gamma}d_{\alpha\delta} = c_\gamma c_\delta, & 1 \leq \gamma, \ \delta \leq s, \ \gamma \neq \delta; \end{cases}$$

*or equivalently $D^T D = \mathrm{Diag}(n_1 + 1, \cdots, n_s + 1) - C^T C$.*

### 4.3. A theorem on residues of meromorphic differential forms.
In this subsection, we review a residue theorem in complex geometry that we will apply in the proof of our main theorems. See [10, Chapter 5] for details.

Let $U$ be the ball $\{z \in \mathbf{C}^s : \|z\| < \varepsilon\}$.

Let $f_1, \cdots, f_s \in \mathcal{O}(\overline{U})$ be functions holomorphic in an open neighborhood of the closure $\overline{U}$ of $U$. We assume that the $f_i$'s have the origin as the isolated common zero.

Set
$$D_i = (f_i) = \text{divisor of } f_i, \quad D = D_1 + \cdots + D_s.$$

Let
$$\omega = \frac{g(z)dz_1 \wedge \cdots \wedge dz_s}{f_1(z) \cdots f_s(z)}$$

be a meromorphic $s$-form with polar divisor $D$.

The **residue** of $\omega$ at the origin is defined as follows,

$$\mathrm{Res}_{(0,\cdots,0)}(\omega) = \left(\frac{1}{2\pi\sqrt{-1}}\right)^s \int_\Gamma \omega,$$

where $\Gamma$ is the real $s$-cycle defined by $\Gamma = \{z : |f_i(z)| = \varepsilon, 1 \leq i \leq s\}$ and is oriented by $d(\arg f_1) \wedge \cdots \wedge d(\arg f_s) \geq 0$.

Let $M$ be a compact complex manifold of dimension $s$.

Suppose that $D_1, \cdots, D_s$ are effective divisors, the intersection of which is a finite set of points.

Let $D = D_1 + \cdots + D_s$.

Let $\omega$ be a meromorphic $s$-form on $M$ with polar divisor $D$.

For each point $P \in D_1 \cap \cdots \cap D_s$, we may restrict $\omega$ to a neighborhood of $U_P$ of $P$ and define the residue $\operatorname{Res}_P(\omega)$ as above.

One has (cf. [10, Chapter 5])

**Lemma 4.1 (Residue Theorem).**
$$\sum_{P \in D_1 \cap \cdots \cap D_s} \operatorname{Res}_P(\omega) = 0$$

**4.4. Proof of Theorem 4.1.** With the above preparations, we are now ready to prove the vanishing theorems.

4.4.1. *Proof of Theorem 4.1, 4k dimensional case.* Let $[V_{(d_{\alpha\beta})}]$ be the fundamental class of $V_{(d_{\alpha\beta})}$ in $H_{4k}(V_{(d_{\alpha\beta})}, \mathbf{Z})$. Then according to (3.1) and the multiplicative property of the Witten genus, up to a constant scalar $\left(\frac{1}{2\pi\sqrt{-1}}\right)^{2k}$,

(4.14)
$$W_c(V_{(d_{\alpha\beta})})$$
$$= \left\langle \left( \left( \prod_{\beta=1}^{s} \left[\frac{i^*x_\beta}{\frac{\theta(i^*x_\beta,\tau)}{\theta'(0,\tau)}}\right]^{n_\beta+1} \right) \left( \prod_{\alpha=1}^{t} \left[\frac{\sum_{\beta=1}^{s} d_{\alpha\beta} i^* x_\beta}{\frac{\theta(\sum_{\beta=1}^{s} d_{\alpha\beta} i^* x_\beta,\tau)}{\theta'(0,\tau)}}\right]^{-1} \right) \frac{\theta_1(i^*c,\tau)}{\theta_1(0,\tau)} \frac{\theta_2(i^*c,\tau)}{\theta_2(0,\tau)} \frac{\theta_3(i^*c,\tau)}{\theta_3(0,\tau)} \right), [V_{(d_{\alpha\beta})}] \right\rangle$$
$$= \left\langle \left( \left( \prod_{\beta=1}^{s} \left[\frac{x_\beta}{\frac{\theta(x_\beta,\tau)}{\theta'(0,\tau)}}\right]^{n_\beta+1} \right) \left( \prod_{\alpha=1}^{t} \left[\frac{1}{\frac{\theta(\sum_{\beta=1}^{s} d_{\alpha\beta} x_\beta,\tau)}{\theta'(0,\tau)}}\right]^{-1} \right) \frac{\theta_1(c,\tau)}{\theta_1(0,\tau)} \frac{\theta_2(c,\tau)}{\theta_2(0,\tau)} \frac{\theta_3(c,\tau)}{\theta_3(0,\tau)} \right), \left[ \prod_{\beta=1}^{s} \mathbf{C}P^{n_\beta} \right] \right\rangle$$
$$= \text{coefficient of } x_1^{n_1} \cdots x_s^{n_s} \text{ in } \left( \frac{(\prod_{\beta=1}^{s} x_\beta^{n_\beta+1}) \left( \prod_{\alpha=1}^{t} \frac{\theta(\sum_{\beta=1}^{s} d_{\alpha\beta} x_\beta,\tau)}{\theta'(0,\tau)} \right)}{\prod_{\beta=1}^{s} \left[\frac{\theta(x_\beta,\tau)}{\theta'(0,\tau)}\right]^{n_\beta+1}} \frac{\theta_1(c,\tau)}{\theta_1(0,\tau)} \frac{\theta_2(c,\tau)}{\theta_2(0,\tau)} \frac{\theta_3(c,\tau)}{\theta_3(0,\tau)} \right)$$
$$= \operatorname{Res}_{(0,\cdots,0)} \left( \frac{\left( \prod_{\alpha=1}^{t} \frac{\theta(\sum_{\beta=1}^{s} d_{\alpha\beta} x_\beta,\tau)}{\theta'(0,\tau)} \right) \frac{\theta_1(c,\tau)}{\theta_1(0,\tau)} \frac{\theta_2(c,\tau)}{\theta_2(0,\tau)} \frac{\theta_3(c,\tau)}{\theta_3(0,\tau)} dx_1 \wedge \cdots \wedge dx_s}{\prod_{\beta=1}^{s} \left[\frac{\theta(x_\beta,\tau)}{\theta'(0,\tau)}\right]^{n_\beta+1}} \right),$$

where we have used Poincaré duality to deduce the second equality.

Set
$$g(x_1,\cdots,x_s) = \left(\prod_{\alpha=1}^{t} \frac{\theta(\sum_{\beta=1}^{s} d_{\alpha\beta}x_\beta,\tau)}{\theta'(0,\tau)}\right) \frac{\theta_1(\sum_{\beta=1}^{s} c_\beta x_\beta,\tau)}{\theta_1(0,\tau)} \frac{\theta_2(\sum_{\beta=1}^{s} c_\beta x_\beta,\tau)}{\theta_2(0,\tau)} \frac{\theta_3(\sum_{\beta=1}^{s} c_\beta x_\beta,\tau)}{\theta_3(0,\tau)},$$

$$f_\beta(x_\beta) = \left[\frac{\theta(x_\beta,\tau)}{\theta'(0,\tau)}\right]^{n_\beta+1}, \quad 1 \leq \beta \leq s,$$

and
$$\omega = \frac{g(x_1,\cdots,x_s)dx_1 \wedge \cdots \wedge dx_s}{f_1(x_1)\cdots f_s(x_s)}.$$

Then up to a constant scalar,

(4.15) $$W_c(V_{(d_{\alpha\beta})}) = \text{Res}_{(0,0,\cdots,0)}(\omega).$$

By (2.14) and (2.16)-(2.18),

(4.16)
$$g(x_1+1, x_2, \cdots, x_s)$$
$$= \left(\prod_{\alpha=1}^{t} \frac{\theta(\sum_{\beta=1}^{s} d_{\alpha\beta}x_\beta + d_{\alpha 1},\tau)}{\theta'(0,\tau)}\right) \frac{\theta_1(\sum_{\beta=1}^{s} c_\beta x_\beta + c_1,\tau)}{\theta_1(0,\tau)} \frac{\theta_2(\sum_{\beta=1}^{s} c_\beta x_\beta + c_1,\tau)}{\theta_2(0,\tau)} \frac{\theta_3(\sum_{\beta=1}^{s} c_\beta x_\beta + c_1,\tau)}{\theta_3(0,\tau)}$$
$$= (-1)^{d_{11}+\cdots+d_{t1}+c_1}\left(\prod_{\alpha=1}^{t} \frac{\theta(\sum_{\beta=1}^{s} d_{\alpha\beta}x_\beta,\tau)}{\theta'(0,\tau)}\right) \frac{\theta_1(\sum_{\beta=1}^{s} c_\beta x_\beta,\tau)}{\theta_1(0,\tau)} \frac{\theta_2(\sum_{\beta=1}^{s} c_\beta x_\beta,\tau)}{\theta_2(0,\tau)} \frac{\theta_3(\sum_{\beta=1}^{s} c_\beta x_\beta,\tau)}{\theta_3(0,\tau)}$$

and
$$f_1(x_1+1) = \left[\frac{\theta(x_1+1,\tau)}{\theta'(0,\tau)}\right]^{n_1+1} = (-1)^{n_1+1}\left[\frac{\theta(x_1,\tau)}{\theta'(0,\tau)}\right]^{n_1+1}.$$

Thus
$$\frac{g(x_1+1,\cdots,x_s)}{f_1(x_1+1)\cdots f_s(x_s)} = (-1)^{(d_{11}+\cdots+d_{t1}+c_1)-(n_1+1)}\frac{g(x_1,\cdots,x_s)}{f_1(x_1)\cdots f_s(x_s)}.$$

Note that by (4.11),
$$(d_{11}+\cdots+d_{t1}+c_1) - (n_1+1) \equiv (d_{11}^2+\cdots+d_{t1}^2+3c_1^2) - (n_1+1) = 0 \bmod 2.$$

Thus one obtains that
$$\frac{g(x_1+1,\cdots,x_s)}{f_1(x_1+1)\cdots f_s(x_s)} = \frac{g(x_1,\cdots,x_s)}{f_1(x_1)\cdots f_s(x_s)}.$$

Similarly, we have

(4.17) $$\frac{g(x_1,\cdots,x_\beta+1,\cdots,x_s)}{f_1(x_1)\cdots f_\beta(x_\beta+1)\cdots f_s(x_s)} = \frac{g(x_1,\cdots,x_s)}{f_1(x_1)\cdots f_s(x_s)}, \quad 1 \leq \beta \leq s.$$

On the other hand, by (2.15) and (2.16)-(2.18),

(4.18)
$$g(x_1 + \tau, x_2, \cdots, x_s)$$
$$= \left( \prod_{\alpha=1}^{t} \frac{\theta(\sum_{\beta=1}^{s} d_{\alpha\beta} x_\beta + d_{\alpha 1}\tau, \tau)}{\theta'(0,\tau)} \right) \frac{\theta_1(\sum_{\beta=1}^{s} c_\beta x_\beta + c_1\tau, \tau)}{\theta_1(0,\tau)} \frac{\theta_2(\sum_{\beta=1}^{s} c_\beta x_\beta + c_1\tau, \tau)}{\theta_2(0,\tau)} \frac{\theta_3(\sum_{\beta=1}^{s} c_\beta x_\beta + c_1\tau, \tau)}{\theta_3(0,\tau)}$$
$$= \left( \prod_{\alpha=1}^{t} (-1)^{d_{\alpha 1}} e^{-2\pi i d_{\alpha 1}(\sum_{\beta=1}^{s} d_{\alpha\beta} x_\beta) - \pi i d_{\alpha 1}^2 \tau} \frac{\theta(\sum_{\beta=1}^{s} d_{\alpha\beta} x_\beta, \tau)}{\theta'(0,\tau)} \right)$$
$$(-1)^{c_1} e^{-6\pi i c_1 (\sum_{\beta=1}^{s} c_\beta x_\beta) - 3\pi i c_1^2 \tau} \frac{\theta_1(\sum_{\beta=1}^{s} c_\beta x_\beta, \tau)}{\theta_1(0,\tau)} \frac{\theta_2(\sum_{\beta=1}^{s} c_\beta x_\beta, \tau)}{\theta_2(0,\tau)} \frac{\theta_3(\sum_{\beta=1}^{s} c_\beta x_\beta, \tau)}{\theta_3(0,\tau)}$$
$$= (-1)^{d_{11} + \cdots + d_{t1} + c_1} e^{-2\pi i \sum_{\alpha=1}^{t} d_{\alpha 1}(\sum_{\beta=1}^{s} d_{\alpha\beta} x_\beta) - \pi i \tau(\sum_{\alpha=1}^{t} d_{\alpha 1}^2) - 6\pi i c_1(\sum_{\beta=1}^{s} c_\beta x_\beta) - 3\pi i c_1^2 \tau}$$
$$\left( \prod_{\alpha=1}^{t} \frac{\theta(\sum_{\beta=1}^{s} d_{\alpha\beta} x_\beta, \tau)}{\theta'(0,\tau)} \right) \frac{\theta_1(\sum_{\beta=1}^{s} c_\beta x_\beta, \tau)}{\theta_1(0,\tau)} \frac{\theta_2(\sum_{\beta=1}^{s} c_\beta x_\beta, \tau)}{\theta_2(0,\tau)} \frac{\theta_3(\sum_{\beta=1}^{s} c_\beta x_\beta, \tau)}{\theta_3(0,\tau)}$$
$$= (-1)^{d_{11} + \cdots + d_{t1} + c_1} e^{-2\pi i \sum_{\alpha=1}^{t} d_{\alpha 1}(\sum_{\beta=1}^{s} d_{\alpha\beta} x_\beta) - \pi i \tau(\sum_{\alpha=1}^{t} d_{\alpha 1}^2) - 6\pi i c_1(\sum_{\beta=1}^{s} c_\beta x_\beta) - 3\pi i c_1^2 \tau} g(x_1, x_2, \cdots, x_s)$$

and

(4.19)
$$f_1(x_1 + \tau) = \left[ \frac{\theta(x_1 + \tau, \tau)}{\theta'(0,\tau)} \right]^{n_1+1}$$
$$= \left[ -e^{-2\pi i x_1 - \pi i \tau} \frac{\theta(x_1, \tau)}{\theta'(0,\tau)} \right]^{n_1+1}$$
$$= (-1)^{n_1+1} e^{-2\pi i (n_1+1) x_1 - \pi i \tau (n_1+1)} \left[ \frac{\theta(x_1, \tau)}{\theta'(0,\tau)} \right]^{n_1+1}$$
$$= (-1)^{n_1+1} e^{-2\pi i (n_1+1) x_1 - \pi i \tau (n_1+1)} f_1(x_1).$$

Therefore
(4.20)
$$\frac{g(x_1 + \tau, \cdots, x_s)}{f_1(x_1 + \tau) \cdots f_s(x_s)}$$
$$= (-1)^{d_{11} + \cdots + d_{t1} + c_1 - n_1 - 1} e^{-2\pi i \sum_{\alpha=1}^{t} d_{\alpha 1}(\sum_{\beta=1}^{s} d_{\alpha\beta} x_\beta) - \pi i \tau(\sum_{\alpha=1}^{t} d_{\alpha 1}^2) - 6\pi i c_1(\sum_{\beta=1}^{s} c_\beta x_\beta) - 3\pi i c_1^2 \tau + 2\pi i (n_1+1) x_1 + \pi i \tau (n_1+1)}$$
$$\cdot \frac{g(x_1, \cdots, x_s)}{f_1(x_1) \cdots f_s(x_s)}.$$

However, we have seen that
$$(d_{11} + \cdots + d_{t1} + c_1) - (n_1 + 1) \equiv (d_{11}^2 + \cdots + d_{t1}^2 + 3c_1^2) - (n_1 + 1) = 0 \bmod 2.$$

Also, one verifies that

$$
\begin{aligned}
&- 2\pi i \sum_{\alpha=1}^{t} d_{\alpha 1}\left(\sum_{\beta=1}^{s} d_{\alpha\beta} x_\beta\right) - \pi i \tau (\sum_{\alpha=1}^{t} d_{\alpha 1}^2) - 6\pi i c_1 (\sum_{\beta=1}^{s} c_\beta x_\beta) \\
&- 3\pi i c_1^2 \tau + 2\pi i (n_1 + 1) x_1 + \pi i \tau (n_1 + 1) \\
(4.21)\quad &= \pi i \tau \left[(n_1+1) - \sum_{\alpha=1}^{t} d_{\alpha 1}^2 - 3 c_1^2\right] + 2\pi i \left[(n_1+1) - \sum_{\alpha=1}^{t} d_{\alpha 1}^2 - 3 c_1^2 \right] x_1 \\
&- 2\pi i \sum_{\beta=2}^{s} \left(\sum_{\alpha=1}^{t} d_{\alpha 1} d_{\alpha\beta} + 3 c_1 c_\beta \right) x_\beta \\
&= 0,
\end{aligned}
$$

where the last equality follows from (4.11).

Consequently, by (4.20), we obtain that

$$
(4.22)\qquad \frac{g(x_1+\tau, \cdots, x_s)}{f_1(x_1+\tau) \cdots f_s(x_s)} = \frac{g(x_1, \cdots, x_s)}{f_1(x_1) \cdots f_s(x_s)}.
$$

Similarly, one also obtains that

$$
(4.23)\qquad \frac{g(x_1, \cdots, x_\beta+\tau, \cdots, x_s)}{f_1(x_1) \cdots f_\beta(x_\beta+\tau) \cdots f_s(x_s)} = \frac{g(x_1, \cdots, x_s)}{f_1(x_1) \cdots f_s(x_s)}, \quad 1 \leq \beta \leq s.
$$

From (4.17) and (4.23), we see that $\omega$ can be viewed as a meromorphic $s$-form defined on the $s$-tori, $(\mathbf{C}/\Gamma)^s$, which is a compact complex manifold.

Recall that $\theta(z, \tau)$ has the lattice points $a+b\tau, a, b \in \mathbf{Z}$ as it's simple zero points (cf. [5]). Thus $\omega$ has pole divisors $\{0\} \times (\mathbf{C}/\Gamma)^{s-1}, (\mathbf{C}/\Gamma) \times \{0\} \times (\mathbf{C}/\Gamma)^{s-2}, \cdots, (\mathbf{C}/\Gamma)^{s-1} \times \{0\}$. It is clear that $(0, 0, \cdots, 0)$ is the unique intersection point of these polar divisors. Therefore by the residue theorem on compact complex manifolds, we directly deduces that

$$
\mathrm{Res}_{(0,0,\cdots,0)}(\omega) = 0.
$$

By (4.15), we obtain that

$$
W_c(V_{(d_{pq})}) = \mathrm{Res}_{(0,0,\cdots,0)}(\omega) = 0.
$$

4.4.2. *Proof of Theorem 4.1, 4k+2 dimensional case.* According to (3.2) and the multiplicative property of the Witten genus, up to a constant scalar $\left(\frac{1}{2\pi\sqrt{-1}}\right)^{2k+1} \sqrt{-1}$,

(4.24)
$W_c(V_{(d_{\alpha\beta})})$

$$= \left\langle \left( \left( \prod_{\beta=1}^{s} \left[ \frac{i^*x_\beta}{\frac{\theta(i^*x_\beta,\tau)}{\theta'(0,\tau)}} \right]^{n_\beta+1} \right) \left( \prod_{\alpha=1}^{t} \left[ \frac{\sum_{\beta=1}^{s} d_{\alpha\beta}i^*x_\beta}{\frac{\theta(\sum_{\beta=1}^{s} d_{\alpha\beta}i^*x_\beta,\tau)}{\theta'(0,\tau)}} \right]^{-1} \right) \frac{\theta(i^*c,\tau)}{\theta_1(0,\tau)\theta_2(0,\tau)\theta_3(0,\tau)} \right), [V_{(d_{\alpha\beta})}] \right\rangle$$

$$= \left\langle \left( \left( \prod_{\beta=1}^{s} \left[ \frac{x_\beta}{\frac{\theta(x_\beta,\tau)}{\theta'(0,\tau)}} \right]^{n_\beta+1} \right) \left( \prod_{\alpha=1}^{t} \left[ \frac{1}{\frac{\theta(\sum_{\beta=1}^{s} d_{\alpha\beta}x_\beta,\tau)}{\theta'(0,\tau)}} \right]^{-1} \right) \frac{\theta(c,\tau)}{\theta_1(0,\tau)\theta_2(0,\tau)\theta_3(0,\tau)} \right), \left[ \prod_{\beta=1}^{s} \mathbf{C}P^{n_\beta} \right] \right\rangle$$

$$= \text{coefficient of } x_1^{n_1}\cdots x_s^{n_s} \text{ in } \left( \frac{(\prod_{\beta=1}^{s} x_\beta^{n_\beta+1}) \left( \prod_{\alpha=1}^{t} \frac{\theta(\sum_{\beta=1}^{s} d_{\alpha\beta}x_\beta,\tau)}{\theta'(0,\tau)} \right)}{\prod_{\beta=1}^{s} \left[ \frac{\theta(x_\beta,\tau)}{\theta'(0,\tau)} \right]^{n_\beta+1}} \frac{\theta(c,\tau)}{\theta_1(0,\tau)\theta_2(0,\tau)\theta_3(0,\tau)} \right)$$

$$= \text{Res}_{(0,\cdots,0)} \left( \frac{\left( \prod_{\alpha=1}^{t} \frac{\theta(\sum_{\beta=1}^{s} d_{\alpha\beta}x_\beta,\tau)}{\theta'(0,\tau)} \right) \frac{\theta(c,\tau)}{\theta_1(0,\tau)\theta_2(0,\tau)\theta_3(0,\tau)} dx_1 \wedge \cdots \wedge dx_s}{\prod_{\beta=1}^{s} \left[ \frac{\theta(x_\beta,\tau)}{\theta'(0,\tau)} \right]^{n_\beta+1}} \right),$$

where we have used the Poincaré duality to deduce the second equality.

Set

$$g(x_1,\cdots,x_s) = \left( \prod_{\alpha=1}^{t} \frac{\theta(\sum_{\beta=1}^{s} d_{\alpha\beta}x_\beta,\tau)}{\theta'(0,\tau)} \right) \frac{\theta(\sum_{\beta=1}^{s} c_\beta x_\beta,\tau)}{\theta_1(0,\tau)\theta_2(0,\tau)\theta_3(0,\tau)},$$

$$f_\beta(x_\beta) = \left[ \frac{\theta(x_\beta,\tau)}{\theta'(0,\tau)} \right]^{n_\beta+1}, \quad 1 \leq \beta \leq s$$

and

$$\omega = \frac{g(x_1,\cdots,x_s)dx_1 \wedge \cdots \wedge dx_s}{f_1(x_1)\cdots f_s(x_s)}.$$

Then up to a constant scalar,

(4.25) $$W_c(V_{(d_{\alpha\beta})}) = \text{Res}_{(0,0,\cdots,0)}(\omega).$$

By (2.14),

(4.26)
$$g(x_1 + 1, x_2, \cdots, x_s) = \left( \prod_{\alpha=1}^{t} \frac{\theta(\sum_{\beta=1}^{s} d_{\alpha\beta} x_\beta + d_{\alpha 1}, \tau)}{\theta'(0, \tau)} \right) \frac{\theta(\sum_{\beta=1}^{s} c_\beta x_\beta + c_1, \tau)}{\theta_1(0, \tau) \theta_2(0, \tau) \theta_3(0, \tau)}$$
$$= (-1)^{d_{11} + \cdots + d_{t1} + c_1} \left( \prod_{\alpha=1}^{t} \frac{\theta(\sum_{\beta=1}^{s} d_{\alpha\beta} x_\beta, \tau)}{\theta'(0, \tau)} \right) \frac{\theta(\sum_{\beta=1}^{s} c_\beta x_\beta, \tau)}{\theta_1(0, \tau) \theta_2(0, \tau) \theta_3(0, \tau)}$$

and

$$f_1(x_1 + 1) = \left[ \frac{\theta(x_1 + 1, \tau)}{\theta'(0, \tau)} \right]^{n_1+1} = (-1)^{n_1+1} \left[ \frac{\theta(x_1, \tau)}{\theta'(0, \tau)} \right]^{n_1+1}.$$

Thus

$$\frac{g(x_1 + 1, \cdots, x_s)}{f_1(x_1 + 1) \cdots f_s(x_s)} = (-1)^{(d_{11} + \cdots + d_{t1} + c_1) - (n_1 + 1)} \frac{g(x_1, \cdots, x_s)}{f_1(x_1) \cdots f_s(x_s)}.$$

Note that by (4.13),

$$(d_{11} + \cdots + d_{t1} + c_1) - (n_1 + 1) \equiv (d_{11}^2 + \cdots + d_{t1}^2 + c_1^2) - (n_1 + 1) = 0 \bmod 2.$$

Thus one obtains that

$$\frac{g(x_1 + 1, \cdots, x_s)}{f_1(x_1 + 1) \cdots f_s(x_s)} = \frac{g(x_1, \cdots, x_s)}{f_1(x_1) \cdots f_s(x_s)}.$$

Similarly, we have

(4.27)
$$\frac{g(x_1, \cdots, x_\beta + 1, \cdots, x_s)}{f_1(x_1) \cdots f_\beta(x_\beta + 1) \cdots f_s(x_s)} = \frac{g(x_1, \cdots, x_s)}{f_1(x_1) \cdots f_s(x_s)}, \quad 1 \leq \beta \leq s.$$

On the other hand, by (2.15),

$$g(x_1 + \tau, x_2, \cdots, x_s)$$

$$= \left( \prod_{\alpha=1}^{t} \frac{\theta(\sum_{\beta=1}^{s} d_{\alpha\beta} x_\beta + d_{\alpha 1}\tau, \tau)}{\theta'(0,\tau)} \right) \frac{\theta(\sum_{\beta=1}^{s} c_\beta x_\beta + c_1\tau, \tau)}{\theta_1(0,\tau)\theta_2(0,\tau)\theta_3(0,\tau)}$$

$$= \left( \prod_{\alpha=1}^{t} (-1)^{d_{\alpha 1}} e^{-2\pi i d_{\alpha 1}(\sum_{\beta=1}^{s} d_{\alpha\beta} x_\beta) - \pi i d_{\alpha 1}^2 \tau} \frac{\theta(\sum_{\beta=1}^{s} d_{\alpha\beta} x_\beta, \tau)}{\theta'(0,\tau)} \right)$$

(4.28)
$$\cdot (-1)^{c_1} e^{-2\pi i c_1(\sum_{\beta=1}^{s} c_\beta x_\beta) - \pi i c_1^2 \tau} \frac{\theta(\sum_{\beta=1}^{s} c_\beta x_\beta, \tau)}{\theta_1(0,\tau)\theta_2(0,\tau)\theta_3(0,\tau)}$$

$$= (-1)^{d_{11}+\cdots+d_{t1}+c_1} e^{-2\pi i \sum_{\alpha=1}^{t} d_{\alpha 1}(\sum_{\beta=1}^{s} d_{\alpha\beta} x_\beta) - \pi i \tau(\sum_{\alpha=1}^{t} d_{\alpha 1}^2) - 2\pi i c_1(\sum_{\beta=1}^{s} c_\beta x_\beta) - \pi i c_1^2 \tau}$$

$$\left( \prod_{\alpha=1}^{t} \frac{\theta(\sum_{\beta=1}^{s} d_{\alpha\beta} x_\beta, \tau)}{\theta'(0,\tau)} \right) \frac{\theta(\sum_{\beta=1}^{s} c_\beta x_\beta, \tau)}{\theta_1(0,\tau)\theta_2(0,\tau)\theta_3(0,\tau)}$$

$$= (-1)^{d_{11}+\cdots+d_{t1}+c_1} e^{-2\pi i \sum_{\alpha=1}^{t} d_{\alpha 1}(\sum_{\beta=1}^{s} d_{\alpha\beta} x_\beta) - \pi i \tau(\sum_{\alpha=1}^{t} d_{\alpha 1}^2) - 2\pi i c_1(\sum_{\beta=1}^{s} c_\beta x_\beta) - \pi i c_1^2 \tau} g(x_1, x_2, \cdots, x_s)$$

and

(4.29)
$$f_1(x_1 + \tau) = \left[ \frac{\theta(x_1 + \tau, \tau)}{\theta'(0,\tau)} \right]^{n_1+1}$$
$$= \left[ -e^{-2\pi i x_1 - \pi i \tau} \frac{\theta(x_1, \tau)}{\theta'(0,\tau)} \right]^{n_1+1}$$
$$= (-1)^{n_1+1} e^{-2\pi i (n_1+1) x_1 - \pi i (n_1+1)\tau} \left[ \frac{\theta(x_1, \tau)}{\theta'(0,\tau)} \right]^{n_1+1}$$
$$= (-1)^{n_1+1} e^{-2\pi i (n_1+1) x_1 - \pi i (n_1+1)\tau} f_1(x_1).$$

Therefore
(4.30)
$$\frac{g(x_1+\tau,\cdots,x_s)}{f_1(x_1+\tau)\cdots f_s(x_s)}$$
$$= (-1)^{d_{11}+\cdots+d_{t1}+c_1-n_1-1} e^{-2\pi i \sum_{\alpha=1}^{t} d_{\alpha 1}(\sum_{\beta=1}^{s} d_{\alpha\beta} x_\beta) - \pi i \tau(\sum_{\alpha=1}^{t} d_{\alpha 1}^2) - 2\pi i c_1(\sum_{\beta=1}^{s} c_\beta x_\beta) - \pi i c_1^2 \tau + 2\pi i (n_1+1) x_1 + \pi i (n_1+1)\tau}$$
$$\cdot \frac{g(x_1,\cdots,x_s)}{f_1(x_1)\cdots f_s(x_s)}.$$

However, we have seen that
$$(d_{11} + \cdots + d_{t1} + c_1) - (n_1 + 1) \equiv (d_{11}^2 + \cdots + d_{t1}^2 + c_1^2) - (n_1 + 1) = 0 \bmod 2.$$

Moreover, one deduces that

$$- 2\pi i \sum_{\alpha=1}^{t} d_{\alpha 1} \left( \sum_{\beta=1}^{s} d_{\alpha\beta} x_\beta \right) - \pi i \tau (\sum_{\alpha=1}^{t} d_{\alpha 1}^2) - 2\pi i c_1 (\sum_{\beta=1}^{s} c_\beta x_\beta)$$
$$- \pi i c_1^2 \tau + 2\pi i (n_1 + 1) x_1 + \pi i \tau (n_1 + 1)$$

(4.31)
$$= \pi i \tau \left[ (n_1 + 1) - \sum_{\alpha=1}^{t} d_{\alpha 1}^2 - c_1^2 \right] + 2\pi i \left[ (n_1 + 1) - \sum_{\alpha=1}^{t} d_{\alpha 1}^2 - c_1^2 \right] x_1$$
$$- 2\pi i \sum_{\beta=2}^{s} \left( \sum_{\alpha=1}^{t} d_{\alpha 1} d_{\alpha\beta} + c_1 c_\beta \right) x_\beta$$
$$= 0,$$

where the last equality follows from (4.13).

Consequently, by (4.30), we obtain that

(4.32) $$\frac{g(x_1 + \tau, \cdots, x_s)}{f_1(x_1 + \tau) \cdots f_s(x_s)} = \frac{g(x_1, \cdots, x_s)}{f_1(x_1) \cdots f_s(x_s)}.$$

Similarly, one also obtains that

(4.33) $$\frac{g(x_1, \cdots, x_\beta + \tau, \cdots, x_s)}{f_1(x_1) \cdots f_\beta(x_\beta + \tau) \cdots f_s(x_s)} = \frac{g(x_1, \cdots, x_s)}{f_1(x_1) \cdots f_s(x_s)}, \quad 1 \leq \beta \leq s.$$

From (4.27) and (4.33), we see that $\omega$ can be viewed as a meromorphic $s$-form defined on the $s$-tori, $(\mathbf{C}/\Gamma)^s$, which is a compact complex manifold.

Recall that $\theta(z, \tau)$ has the lattice points $a + b\tau, a, b \in \mathbf{Z}$ as it's simple zero points [5]. We see then that $\omega$ has pole divisors $\{0\} \times (\mathbf{C}/\Gamma)^{s-1}, (\mathbf{C}/\Gamma) \times \{0\} \times (\mathbf{C}/\Gamma)^{s-2}, \cdots, (\mathbf{C}/\Gamma)^{s-1} \times \{0\}$. Thus $(0, 0, \cdots, 0)$ is the unique intersection point of these polar divisors.

By the residue theorem on compact complex manifolds, we directly deduce that

$$\mathrm{Res}_{(0,0,\cdots,0)}(\omega) = 0.$$

By (4.25), we obtain the desired result that

$$W_c(V_{(d_{\alpha\beta})}) = \mathrm{Res}_{(0,0,\cdots,0)}(\omega) = 0.$$

### 4.5. Proof of Theorem 4.2 and Corollary 4.2.
To prove Theorem 4.2, we will need the Rokhlin type congruence established in [24, 27].

Let $M$ be an $8k+4$ dimensional spin$^c$ manifold. Let $B \subset M$ be an $8k+2$ dimensional closed oriented submanifold of $M$ such that $[B] \in H_{8k+2}(M, \mathbf{Z_2})$ is dual to the Stiefel-Whitney class $w_2(TM)$. In such a situation $B$ is called a characteristic submanifold of $M$. The existence of $B$ is clear. Moreover, $M \setminus B$ is spin. We fix a spin structure on $M \setminus B$. Then $B$ carries a canonically induced spin structure. Furthermore, the spin cobordism class of the induced spin structure does not depend on the spin structure chosen on $M \setminus B$.

Let $i_B : B \hookrightarrow M$ denote the canonical embedding of $B$ into $M$.

Let $e \in H^2(M, \mathbf{Z})$ be the dual of $[B] \subset H_{8k+2}(M, \mathbf{Z})$.

Let $E$ be a real vector bundle over $M$. Then $i_B^* E$ is a real vector bundle over the spin manifold $B$. Let $\mathrm{ind}_2(i_B^* E)$ be the mod 2 index in the sense of Atiyah-Singer [4] associated to $i_B^* E$. Then we have the following analytic Rokhlin congruence.

**Theorem 4.3** (Zhang [24], [27]). *The following identity holds,*

(4.34) $$\left\langle \widehat{A}(TM) \mathrm{ch}(E \otimes \mathbf{C}) \exp\left(\frac{e}{2}\right), [M] \right\rangle \equiv \mathrm{ind}_2(i_B^* E) \mod 2.$$

With the above preparation, we now prove Theorem 4.2 in the following.

Since $V_{(d_{\alpha\beta})}$ is string, one has

(4.35) $$p_1(V_{(d_{\alpha\beta})}) = 0.$$

By proceeding similarly as in Propositions 4.1 and 4.2, one easily gets the following identities:

(4.36) $$\begin{cases} n_\beta + 1 - \sum_{\alpha=1}^{t} d_{\alpha\beta}^2 = 0, & 1 \leq \beta \leq s, \\ \sum_{\alpha=1}^{t} d_{\alpha\gamma} d_{\alpha\delta} = 0, & 1 \leq \gamma, \delta \leq s, \ \gamma \neq \delta. \end{cases}$$

Since in $D$ there is a row consisting of only even numbers, without loss of generality, assume $d_{t\beta}, 1 \leq \beta \leq s$, are even numbers.

Let $V_t$ be the generalized complete intersection determined by $\overset{s}{\underset{\beta=1}{\otimes}} P_\beta^*(\mathcal{O}_\beta(d_{\alpha\beta}))$, $1 \leq \alpha \leq t-1$, which is an $8k+4$ dimensional generalized complete intersection.

Denote $\overset{s}{\underset{\beta=1}{\otimes}} P_\beta^*(\mathcal{O}_\beta(d_{t\beta}))$ by $\zeta$.

Let $i : V_t \hookrightarrow \mathbf{C}P^{n_1} \times \cdots \times \mathbf{C}P^{n_s}$ and $i_{V_{(d_{\alpha\beta})}} : V_{(d_{\alpha\beta})} \hookrightarrow V_t$ be the notations of embeddings.

Clearly, $[V_{(d_{\alpha\beta})}] \in H_{8k+2}(V_t, \mathbf{Z}_2)$ is dual to $c_1(i^*\zeta)$ mod 2. However,

$$c_1(i^*\zeta) = \sum_{\beta=1}^{s} d_{t\beta} i^* x_\beta \equiv \sum_{\beta=1}^{s} d_{t\beta}^2 i^* x_\beta \mod 2$$

and by (4.36),

$$\sum_{\beta=1}^{s} d_{t\beta}^2 i^* x_\beta = \sum_{\beta=1}^{s} (n_\beta + 1 - \sum_{\alpha=1}^{t-1} d_{\alpha\beta}^2) i^* x_\beta \equiv \sum_{\beta=1}^{s} (n_\beta + 1 - \sum_{\alpha=1}^{t-1} d_{\alpha\beta}) i^* x_\beta \equiv w_2(TV_t) \mod 2.$$

So $[V_{(d_{\alpha\beta})}] \in H_{8k+2}(V_t, \mathbf{Z}_2)$ is dual to the second Stiefel-Whitney class of $V_t$, i.e., $V_{(d_{\alpha\beta})}$ is a characteristic submanifold of $V_t$.

One can then apply Zhang's Rokhlin type congruence to get

(4.37) $$\phi(V_{(d_{\alpha\beta})}) = \mathrm{ind}_2(\Theta(TV_{(d_{\alpha\beta})})) \\ \equiv \left\langle \widehat{A}(TV_t) \cosh\left(\frac{\sum_{\beta=1}^{s} d_{t\beta} i^* x_\beta}{2}\right) \mathrm{ch}(\Theta(T_\mathbf{C}V_t - i^*\zeta_\mathbf{R} \otimes \mathbf{C})), [V_t] \right\rangle \mod 2\mathbf{Z}[[q]],$$

where we have used $i^*_{V_{(d_{\alpha\beta})}}(TV_t - i^*\zeta_\mathbf{R}) = TV_{(d_{\alpha\beta})}$.

To continue our calculation, we need the following lemma:

**Lemma 4.2.** *The following identity holds,*

(4.38) $$\Theta(T_\mathbf{C}V_t - i^*\zeta_\mathbf{R} \otimes \mathbf{C}) \equiv \Theta(T_\mathbf{C}V_t) \otimes \left( \overset{\infty}{\underset{n=1}{\otimes}} \Lambda_{q^{2n}}(\widetilde{i^*\zeta_\mathbf{R} \otimes \mathbf{C}}) \right) \mod 2(\widetilde{i^*\zeta_\mathbf{R} \otimes \mathbf{C}}) \cdot K(V_t)[[q]].$$

*Proof.* By using a trick as in [12, (3.10)], one deduces that

$$\begin{aligned}
\Theta(T_{\mathbf{C}}V_t - i^*\zeta_{\mathbf{R}} \otimes \mathbf{C}) &= \overset{\infty}{\underset{m=1}{\otimes}} S_{q^{2m}}\left(\widetilde{T_{\mathbf{C}}V_t} - i^*\widetilde{\zeta_{\mathbf{R}} \otimes \mathbf{C}}\right) \\
&= \overset{\infty}{\underset{m=1}{\otimes}} \left(\frac{1}{\Lambda_{-q^{2m}}(\widetilde{T_{\mathbf{C}}V_t} - i^*\widetilde{\zeta_{\mathbf{R}} \otimes \mathbf{C}})}\right) \\
&= \overset{\infty}{\underset{m=1}{\otimes}} \left(\frac{\Lambda_{-q^{2m}}(i^*\widetilde{\zeta_{\mathbf{R}} \otimes \mathbf{C}})}{\Lambda_{-q^{2m}}(\widetilde{T_{\mathbf{C}}V_t})}\right) \\
&\equiv \overset{\infty}{\underset{m=1}{\otimes}} \left(\frac{\Lambda_{q^{2m}}(i^*\widetilde{\zeta_{\mathbf{R}} \otimes \mathbf{C}})}{\Lambda_{-q^{2m}}(\widetilde{T_{\mathbf{C}}V_t})}\right) \\
&= \Theta(T_{\mathbf{C}}V_t) \otimes \left(\overset{\infty}{\underset{n=1}{\otimes}} \Lambda_{q^{2n}}(i^*\widetilde{\zeta_{\mathbf{R}} \otimes \mathbf{C}})\right) \mod 2(i^*\widetilde{\zeta_{\mathbf{R}} \otimes \mathbf{C}}) \cdot \mathbf{Z}[T_{\mathbf{C}}V_t, i^*\zeta_{\mathbf{R}} \otimes \mathbf{C}][[q]].
\end{aligned}$$

□

From (4.37) and (4.38), we have

(4.39)
$$\phi(V_{(d_{\alpha\beta})})$$
$$\equiv \left\langle \widehat{A}(TV_t) \cosh\left(\frac{\sum_{\beta=1}^{s} d_{t\beta} i^* x_\beta}{2}\right) \mathrm{ch}\left(\bigotimes_{n=1}^{\infty} \Lambda_{q^{2n}}(i^*\widetilde{\zeta_{\mathbf{R}} \otimes \mathbf{C}})\right) \mathrm{ch}\left(\Theta(T_{\mathbf{C}} V_t)\right), [V_t]\right\rangle$$

$$= \left\langle \widehat{A}(TV_t)\mathrm{ch}(\Theta(T_{\mathbf{C}} V_t)) \cdot \frac{\theta_1\left(\frac{\sum_{\beta=1}^{s} d_{t\beta} i^* x_\beta}{2\pi i}, \tau\right)}{\theta_1(0,\tau)}, [V_t]\right\rangle$$

$$= \left\langle \prod_{\beta=1}^{s}\left[\frac{x_\beta}{2\pi i \frac{\theta\left(\frac{x_\beta}{2\pi i},\tau\right)}{\theta'(0,\tau)}}\right]^{n_\beta+1} \left(\prod_{\alpha=1}^{t-1} \frac{2\pi i\theta\left(\frac{\sum_{\beta=1}^{s} d_{\alpha\beta} x_\beta}{2\pi i},\tau\right)}{\theta'(0,\tau)}\right) \frac{\theta_1\left(\frac{\sum_{\beta=1}^{s} d_{t\beta} x_\beta}{2\pi i},\tau\right)}{\theta_1(0,\tau)}, \left[\prod_{\beta=1}^{s} \mathbf{C}P^{n_\beta}\right]\right\rangle$$

$$= \mathrm{Res}_{(0,\cdots,0)}\left(\prod_{\beta=1}^{s}\left[\frac{\theta'(0,\tau)}{2\pi i\theta(\frac{x_\beta}{2\pi i},\tau)}\right]^{n_\beta+1} \left(\prod_{\alpha=1}^{t-1} \frac{2\pi i\theta\left(\frac{\sum_{\beta=1}^{s} d_{\alpha\beta} x_\beta}{2\pi i},\tau\right)}{\theta'(0,\tau)}\right)\right.$$
$$\left. \cdot \frac{\theta_1\left(\frac{\sum_{\beta=1}^{s} d_{t\beta} x_\beta}{2\pi i},\tau\right)}{\theta_1(0,\tau)} dx_1 \wedge dx_2 \wedge \ldots \wedge dx_s\right) \mod 2\mathbf{Z}[[q]].$$

Set
$$g(x_1,\cdots,x_s) = \prod_{\alpha=1}^{t-1}\left(\frac{2\pi i\theta\left(\frac{\sum_{\beta=1}^{s} d_{\alpha\beta} x_\beta}{2\pi i},\tau\right)}{\theta'(0,\tau)}\right) \frac{\theta_1\left(\frac{\sum_{\beta=1}^{s} d_{t\beta} x_\beta}{2\pi i},\tau\right)}{\theta_1(0,\tau)},$$

$$f_\beta(x_\beta) = \left[\frac{2\pi i\theta(\frac{x_\beta}{2\pi i},\tau)}{\theta'(0,\tau)}\right]^{n_\beta+1}, \quad 1 \leq \beta \leq s,$$

and
$$\omega = \frac{g(x_1,\cdots,x_s)dx_1 \wedge \cdots \wedge dx_s}{f_1(x_1)\cdots f_s(x_s)}.$$

Then we have,

(4.40) $$\phi(V_{(d_{\alpha\beta})}) \equiv \mathrm{Res}_{(0,0,\cdots,0)}(\omega) \mod 2.$$

By (2.14) and (2.16)-(2.18),

(4.41)
$$g(x_1 + 2\pi i, x_2, \cdots, x_s) = \left( \prod_{\alpha=1}^{t-1} \frac{2\pi i \theta\left(\frac{\sum_{\beta=1}^{s} d_{\alpha\beta}x_\beta}{2\pi i} + d_{\alpha 1}, \tau\right)}{\theta'(0,\tau)} \right) \cdot \frac{\theta_1\left(\frac{\sum_{\beta=1}^{s} d_{t\beta}x_\beta}{2\pi i} + d_{t1}, \tau\right)}{\theta_1(0,\tau)}$$

$$= (-1)^{d_{11}+\cdots+d_{t1}} \left( \prod_{\alpha=1}^{t-1} \frac{2\pi i \theta\left(\frac{\sum_{\beta=1}^{s} d_{\alpha\beta}x_\beta}{2\pi i}, \tau\right)}{\theta'(0,\tau)} \right) \frac{\theta_1\left(\frac{\sum_{\beta=1}^{s} d_{t\beta}x_\beta}{2\pi i}, \tau\right)}{\theta_1(0,\tau)}$$

and

$$f_1(x_1 + 2\pi i) = \left[\frac{2\pi i \theta(\frac{x_1}{2\pi i} + 1, \tau)}{\theta'(0,\tau)}\right]^{n_1+1} = (-1)^{n_1+1}\left[\frac{2\pi i \theta(\frac{x_1}{2\pi i}, \tau)}{\theta'(0,\tau)}\right]^{n_1+1}.$$

Thus

$$\frac{g(x_1 + 2\pi i, \cdots, x_s)}{f_1(x_1 + 2\pi i) \cdots f_s(x_s)} = (-1)^{(d_{11}+\cdots+d_{t1})-(n_1+1)} \frac{g(x_1, \cdots, x_s)}{f_1(x_1) \cdots f_s(x_s)}.$$

Note that by (4.36),

$$(d_{11} + \cdots + d_{t1}) - (n_1 + 1) \equiv (d_{11}^2 + \cdots + d_{t1}^2) - (n_1 + 1) = 0 \mod 2.$$

Thus one obtains that

$$\frac{g(x_1 + 2\pi i, \cdots, x_s)}{f_1(x_1 + 2\pi i) \cdots f_s(x_s)} = \frac{g(x_1, \cdots, x_s)}{f_1(x_1) \cdots f_s(x_s)}.$$

Similarly, we have

(4.42) $$\frac{g(x_1, \cdots, x_\beta + 2\pi i, \cdots, x_s)}{f_1(x_1) \cdots f_\beta(x_\beta + 2\pi i) \cdots f_s(x_s)} = \frac{g(x_1, \cdots, x_s)}{f_1(x_1) \cdots f_s(x_s)}, \quad 1 \leq \beta \leq s.$$

On the other hand, by (2.15) and (2.16)-(2.18),

(4.43)
$$g(x_1 + 2\pi i\tau, x_2, \cdots, x_s) = \left(\prod_{\alpha=1}^{t-1} \frac{2\pi i\theta\left(\frac{\sum_{\beta=1}^{s} d_{\alpha\beta}x_\beta}{2\pi i} + d_{\alpha 1}\tau, \tau\right)}{\theta'(0,\tau)}\right) \frac{\theta_1\left(\frac{\sum_{\beta=1}^{s} d_{t\beta}x_\beta}{2\pi i} + d_{t 1}\tau, \tau\right)}{\theta_1(0,\tau)}$$

$$= \left(\prod_{\alpha=1}^{t-1} (-1)^{d_{\alpha 1}} e^{-d_{\alpha 1}\left(\sum_{\beta=1}^{s} d_{\alpha\beta}x_\beta\right) - \pi i d_{\alpha 1}^2 \tau} \frac{2\pi i\theta\left(\frac{\sum_{\beta=1}^{s} d_{\alpha\beta}x_\beta}{2\pi i}, \tau\right)}{\theta'(0,\tau)}\right)$$

$$\cdot e^{-d_{t 1}\sum_{\beta=1}^{s} d_{t\beta}x_\beta - \pi i d_{t 1}^2 \tau} \frac{\theta_1\left(\frac{\sum_{\beta=1}^{s} d_{t\beta}x_\beta}{2\pi i}, \tau\right)}{\theta_1(0,\tau)}$$

$$= (-1)^{d_{11}+\cdots+d_{(t-1)1}} e^{-\sum_{\alpha=1}^{t-1} d_{\alpha 1}\left(\sum_{\beta=1}^{s} d_{\alpha\beta}x_\beta\right) - d_{t 1}\sum_{\beta=1}^{s} d_{t\beta}x_\beta - \pi i\left(\sum_{\alpha=1}^{t-1} d_{\alpha 1}^2\right)\tau - \pi i d_{t 1}^2 \tau}$$

$$\cdot \left(\prod_{\alpha=1}^{t-1} \frac{2\pi i\theta\left(\frac{\sum_{\beta=1}^{s} d_{\alpha\beta}x_\beta}{2\pi i}, \tau\right)}{\theta'(0,\tau)}\right) \cdot \frac{\theta_1\left(\frac{\sum_{\beta=1}^{s} d_{t\beta}x_\beta}{2\pi i}, \tau\right)}{\theta_1(0,\tau)}$$

$$= (-1)^{d_{11}+\cdots+d_{(t-1)1}} e^{-\sum_{\alpha=1}^{t} d_{\alpha 1}\left(\sum_{\beta=1}^{s} d_{\alpha\beta}x_\beta\right) - \pi i\left(\sum_{\alpha=1}^{t} d_{\alpha 1}^2\right)\tau} g(x_1, x_2, \cdots, x_s)$$

and

(4.44)
$$f_1(x_1 + 2\pi i\tau) = \left[\frac{2\pi i\theta(\frac{x_1}{2\pi i} + \tau, \tau)}{\theta'(0,\tau)}\right]^{n_1+1} = \left[-e^{-x_1 - \pi i\tau}\frac{2\pi i\theta(\frac{x_1}{2\pi i}, \tau)}{\theta'(0,\tau)}\right]^{n_1+1}$$

$$= (-1)^{n_1+1} e^{-(n_1+1)x_1 - \pi i(n_1+1)\tau}\left[\frac{2\pi i\theta(\frac{x_1}{2\pi i}, \tau)}{\theta'(0,\tau)}\right]^{n_1+1}$$

$$= (-1)^{n_1+1} e^{-(n_1+1)x_1 - \pi i(n_1+1)\tau} f_1(x_1).$$

Therefore

(4.45)
$$\frac{g(x_1 + 2\pi i\tau, \cdots, x_s)}{f_1(x_1 + 2\pi i\tau)\cdots f_s(x_s)}$$

$$= (-1)^{d_{11}+\cdots+d_{(t-1)1}-n_1-1} e^{-\sum_{\alpha=1}^{t} d_{\alpha 1}\left(\sum_{\beta=1}^{s} d_{\alpha\beta}x_\beta\right) - \pi i\left(\sum_{\alpha=1}^{t} d_{\alpha 1}^2\right)\tau + (n_1+1)x_1 + \pi i(n_1+1)\tau}$$

$$\cdot \frac{g(x_1, \cdots, x_s)}{f_1(x_1)\cdots f_s(x_s)}.$$

Now by (4.36) and the assumption that $d_{t\beta}, 1 \leq \beta \leq s$, are even numbers, one has

$$(d_{11} + \cdots + d_{(t-1)1}) - (n_1 + 1) \equiv (d_{11}^2 + \cdots + d_{t1}^2) - (n_1 + 1) = 0 \mod 2$$

and

$$
\begin{aligned}
&-\sum_{\alpha=1}^{t} d_{\alpha 1} \left( \sum_{\beta=1}^{s} d_{\alpha\beta} x_\beta \right) - \pi i (\sum_{\alpha=1}^{t} d_{\alpha 1}^2) \tau + (n_1 + 1) x_1 + \pi i (n_1 + 1) \tau \\
&= \pi i \left[ (n_1 + 1) - \sum_{\alpha=1}^{t} d_{\alpha 1}^2 \right] \tau + \left[ (n_1 + 1) - \sum_{\alpha=1}^{t} d_{\alpha 1}^2 \right] x_1 - \sum_{\beta=2}^{s} \left( \sum_{\alpha=1}^{t} d_{\alpha 1} d_{\alpha\beta} \right) x_\beta \\
&= 0.
\end{aligned}
\tag{4.46}
$$

Consequently, by (4.45), we obtain that

$$\frac{g(x_1 + 2\pi i \tau, \cdots, x_s)}{f_1(x_1 + 2\pi i \tau) \cdots f_s(x_s)} = \frac{g(x_1, \cdots, x_s)}{f_1(x_1) \cdots f_s(x_s)}. \tag{4.47}$$

Similarly, one also obtains that

$$\frac{g(x_1, \cdots, x_\beta + 2\pi i \tau, \cdots, x_s)}{f_1(x_1) \cdots f_\beta(x_\beta + 2\pi i \tau) \cdots f_s(x_s)} = \frac{g(x_1, \cdots, x_s)}{f_1(x_1) \cdots f_s(x_s)}, \quad 1 \leq \beta \leq s. \tag{4.48}$$

Therefore from (4.42) and (4.48), we see that $\omega$ can be viewed as a meromorphic $s$-form defined on the $s$-tori, $(\mathbf{C}/\Gamma)^s$, which is a compact complex manifold. Here $\Gamma = \{(a + b\tau)2\pi i\}$.

We know $(0, 0, \cdots, 0)$ is its only pole. Therefore by the residue theorem on compact complex manifolds, we directly deduce that

$$\mathrm{Res}_{(0,0,\cdots,0)}(\omega) = 0.$$

By (4.40), we obtain the desired result that

$$\phi(V_{(d_{\alpha\beta})}) \equiv \mathrm{Res}_{(0,0,\cdots,0)}(\omega) = 0 \mod 2.$$

At last, to prove Corollary 4.2, observe that if $V(d_1, \cdots, d_r) \subset \mathbf{C}P^{4k+1+r}$ is a string complete intersection of complex dimension $4k + 1$ ($k \geq 1$) in a complex projective space, then we have (by (4.36)),

$$4k + 1 + r + 1 - \sum_{\alpha=1}^{r} d_\alpha^2 = 0$$

or

$$\sum_{\alpha=1}^{r} (d_\alpha^2 - 1) = 4k + 2. \tag{4.49}$$

If each $d_\alpha, 1 \leq \alpha \leq r$, is an odd number, then the left-hand side of (4.49) is divisible by 8, which results in a contradiction. So at least one of $d_\alpha$'s, $1 \leq \alpha \leq r$, has to be even. Then Corollary 4.2 follows directly from Theorem 4.2.

**Acknowledgement** We are indebted to Kefeng Liu for very helpful discussions. The work of the second author was partially supported by a start-up grant from National University of Singapore. The work of the third author was partially supported by NNSFC and MOEC.

Q. Chen, Department of Mathematics, University of Southern California, Los Angeles, CA 90089, USA (qingtaoc@usc.edu)

F. Han, Department of Mathematics, National University of Singapore, Block S17, 10 Lower Kent Ridge Road, Singapore 119076 (mathanf@nus.edu.sg)

W. Zhang, Chern Institute of Mathematics & LPMC, Nankai University, Tianjin 300071, P.R. China. (weiping@nankai.edu.cn)